\documentclass[11pt]{amsart}
\usepackage{amssymb}
\usepackage{amsfonts}
\usepackage{color}
\usepackage[pdftex]{graphicx}

\newcommand{\nc}{\newcommand}
\nc{\thref}[1]{Theorem~\ref{theo:#1}}
\nc{\selabel}[1]{\label{sect:#1}}
\nc{\seref}[1]{Section~\ref{sect:#1}}
\nc{\lelabel}[1]{\label{lemm:#1}}
\nc{\leref}[1]{Lemma~\ref{lemm:#1}}
\nc{\prlabel}[1]{\label{prop:#1}}
\nc{\prref}[1]{Proposition~\ref{prop:#1}}
\nc{\colabel}[1]{\label{coro:#1}}
\nc{\coref}[1]{Corollary~\ref{coro:#1}}
\nc{\exlabel}[1]{\label{exam:#1}}
\nc{\exref}[1]{Example~\ref{exam:#1}}
\nc{\delabel}[1]{\label{defi:#1}}
\nc{\deref}[1]{Definition~\ref{defi:#1}}
\nc{\eqlabel}[1]{\label{equa:#1}}
\nc{\relabel}[1]{\label{rema:#1}}
\nc{\reref}[1]{Lemma~\ref{rema:#1}}
\providecommand{\operatorname}[1]{\mathrm{#1}\,}
\nc{\Hom}{\operatorname{Hom}}
\nc{\Mor}{\operatorname{Mor}}
\nc{\Aut}{\operatorname{Aut}}
\nc{\Ann}{\operatorname{Ann}}
\nc{\Ker}{\operatorname{Ker}}
\nc{\Trace}{\operatorname{Trace}}
\nc{\Char}{\operatorname{Char}}
\nc{\Mod}{\operatorname{Mod}}
\nc{\End}{\operatorname{End}}
\nc{\Spec}{\operatorname{Spec}}
\nc{\Span}{\operatorname{Span}}
\nc{\sgn}{\operatorname{sgn}}
\nc{\Id}{\operatorname{Id}}
\nc{\Com}{\operatorname{Com}}
\nc{\rank}{\operatorname{rank}}

\let\:=\colon

\newtheorem{de}{Definition}[section]
\newtheorem{lm}[de]{Lemma}
\newtheorem{pr}[de]{Proposition}
\newtheorem{co}[de]{Corollary}
\newtheorem{re}[de]{Remark}
\newtheorem{res}[de]{Remarks}
\newtheorem{te}[de]{Theorem}
\newtheorem{ex}[de]{Example}
\newtheorem{exs}[de]{Examples}

\def\bex{\begin{ex}}
\def\eex{\end{ex}}
\def\bexs{\begin{exs}}
\def\eexs{\end{exs}}
\def\bl{\begin{lm}}
\def\el{\end{lm}}
\def\bc{\begin{co}}
\def\ec{\end{co}}
\def\bt{\begin{te}}
\def\et{\end{te}}
\def\bpr{\begin{pr}}
\def\epr{\end{pr}}
\def\br{\begin{re}}
\def\er{\end{re}}
\def\brs{\begin{res}}
\def\ers{\end{res}}
\def\bd{\begin{de}}
\def\ed{\end{de}}
\def\be{\begin{equation}}
\def\ee{\end{equation}}
\def\bea{\begin{eqnarray*}}
\def\eea{\end{eqnarray*}}
\def\bp{\begin{proof}}
\def\ep{\end{proof}}

\let\:=\colon
\begin{document}
\title{Generalization of the Apollonius Circles}
\author{Cosmin Pohoata and Vladimir Zajic \\}

\date{}
\maketitle
\small
\begin{abstract}
The three Apollonius circles of a triangle, each passing through a triangle vertex, the corresponding vertex of the cevian triangle of the incenter and the corresponding vertex of the circumcevian triangle of the symmedian point, are coaxal. Similarly defined three circles remain coaxal, when the circumcevian triangle is defined with respect to any point on the triangle circumconic through the incenter and symmedian point. Inversion in the incircle of the reference triangle carries these three coaxal circles into coaxal circles, each passing through a vertex of the inverted triangle and centered on the opposite sideline, at the intersection of the orthotransversal with respect to a point on the Euler line of the inverted triangle. A similar circumconic exists in a more general configuration, when the cevian triangle is defined with respect to an arbitrary point, passing through this arbitrary point and isogonal conjugate of its complement.
\end{abstract}
\bigskip

\bigskip

{\bf{1. Preliminaries}}

\bigskip

Let $X$, $Y$, $Z$ be the feet of the internal bisectors and $X^{\prime}$, $Y^{\prime}$, $Z^{\prime}$ the feet of the external bisectors of the angles $\angle A$, $\angle B$, $\angle C$ in a triangle $ABC$. The three Apollonius circles $(AXX^{\prime})$, $(BYY^{\prime})$, $(CZZ^{\prime})$ of the triangle, each passing through a triangle vertex and its corresponding feet of the internal and external angle bisectors, are loci of points with equal ratios of distances from the remaining two triangle vertices. By the transitivity of equivalence, they are coaxal with two common points, the two isodynamic points of the triangle. Moreover, the radical axes of the three Apollonius circles with the triangle circumcircle coincide with the corresponding symmedians of the triangle [1].

Let $A^{\prime}B^{\prime}C^{\prime}$ be the circumcevian triangle of the reference triangle with respect to a point $Q$. If $Q = K$ is the symmedian point, the Apollonius circles $(AXA^{\prime})$, $(BYB^{\prime})$, $(CZC^{\prime})$ are coaxal. In this paper, we study the locus of points $Q$, for which the circles $(AXA^{\prime})$, $(BYB^{\prime})$, $(CZC^{\prime})$ remain coaxal. At first, we synthetically prove the coaxality for four special cases, and afterwards, we make use of the homogenous barycentric coordinates to prove a more general proposition, communicated by Jean-Pierre Ehrmann, in which the feet of the internal angle bisectors are replaced with the traces $X$, $Y$, $Z$ of an arbitrary point $P$. In order to set up our main theorem, we begin with four simple lemmas:

\bigskip

{\bf{Lemma 1}}. Let $P$ be an arbitrary point in the plane of a triangle $ABC$ and $XYZ$ the cevian triangle of the point $P$. Let three arbitrary circles $\mathcal U$, $\mathcal V$, $\mathcal W$, passing through the endpoints of the segments $AX$, $BY$, $CZ$, meet the triangle circumcircle $\mathcal O$ again at points $A^{\prime}$, $B^{\prime}$, $C^{\prime}$ and the triangle sidelines $BC$, $CA$, $AB$ again at points $U$, $V$, $W$, respectively. The points $U$, $V$, $W$ are collinear, if and only if the lines $AA^{\prime}$, $BB^{\prime}$, $CC^{\prime}$ are concurrent. 

\bigskip

{\it{Proof of Lemma 1}}. Let the circumcircle secants $AA^{\prime}$, $ BB^{\prime}$, $CC^{\prime}$ cut the triangle sidelines $ BC$, $CA$, $AB$ at points $Q_{a}$, $Q_{b}$, $Q_{c}$. The line $AA^{\prime}$ is the radical axis of the circumcircle $\mathcal O$ and circle $\mathcal U$. Powers of the point $Q_{a}$ to both circles are equal,

$$\overline{Q_{a}C} \cdot \overline{Q_{a}B} = \overline{Q_{a}U} \cdot \overline{Q_{a}X}.$$

Using this and $\overline{XC} - \overline{XB} = \overline{UC} - \overline{UB} = \overline{Q_{a}C} - \overline{Q_{a}B} = \overline{BC}$ leads to

$$\frac{\overline{UB}}{\overline{UC}} = \frac{\overline{Q_{a}B}}{\overline{Q_{a}C}} \cdot \frac{\overline{XC}}{\overline{XB}}.$$

Multiplying the results, after cyclic exchange,

$$\frac{\overline{UB}}{\overline{UC}} \cdot \frac{\overline{VC}}{\overline{VA}} \cdot \frac{\overline{WA}}{\overline{WB}} = \left(\frac{\overline{Q_{a}B}}{\overline{Q_{a}C}} \cdot \frac{\overline{XC}}{\overline{XB}}\right) \cdot \left(\frac{\overline{Q_{b}C}}{\overline{Q_{b}A}} \cdot \frac{\overline{YA}}{\overline{YC}}\right) \cdot \left(\frac{\overline{Q_{c}A}}{\overline{Q_{c}B}} \cdot \frac{\overline{ZB}}{\overline{ZA}}\right).$$

By Ceva and Menelaus theorems, the points $U,$ $V,$ $W$ are collinear, if and only if the lines $AA^{\prime}$, $BB^{\prime}$, $CC^{\prime}$ are concurrent at a point $Q$. $\square$

\bigskip

{\bf{Lemma 2}}. Let $A_{0}B_{0}C_{0}$ be a triangle centrally similar to a reference triangle $ABC$, with similarity center at the incenter $I$ of the reference triangle. Let the external bisectors of the angles $\angle A_{0}$, $\angle B_{0}$, $\angle C_{0}$ cut the sidelines $BC$, $CA,$ $AB$ of the reference triangle at points $U$, $V$, $W$. The points $U$, $V$, $W$ are collinear.

\bigskip

{\it{Proof of Lemma 2}}. The internal angle bisectors of the triangles $A_{0}B_{0}C_{0}$, $ABC$ are identical. Pairs of the external angle bisectors $A_{0}U$, $B_{0}V$, $C_{0}W$ of the triangle $A_{0}B_{0}C_{0}$ meet on the corresponding internal angle bisectors at the three excenters $I_{a}$, $I_{b}$, $I_{c}$ of this triangle. In addition, the external angle bisectors $A_{0}U$, $B_{0}V$, $C_{0}W$ cut the opposite sidelines of the triangle $A_{0}B_{0}C_{0}$ at points $X_{0}^{\prime}$, $Y_{0}^{\prime}$, $Z_{0}^{\prime}$, collinear by Menelaus theorem. The triangles $CUV$, $C_{0}X_{0}^{\prime}Y_{0}^{\prime}$ are centrally similar with similarity center $I_{c}$, having the corresponding sidelines $CU$, $C_{0}X_{0}^{\prime}$ and $CV$, $C_{0}Y_{0}^{\prime}$ parallel and the lines $CC_{0}$, $UX_{0}^{\prime}$, $VY_{0}^{\prime}$, connecting their corresponding vertices, concurrent at their similarity center $I_{c}$. It follows that the lines $UV$, $X_{0}^{\prime}Y_{0}^{\prime}$ are parallel and similarly, the lines $VW$, $Y_{0}^{\prime}Z_{0}^{\prime}$ are also parallel. Since the points $X_{0}^{\prime}$, $Y_{0}^{\prime}$, $Z_{0}^{\prime}$ are collinear, the points $U$, $V$, $W$ are also collinear. $\square$

\bigskip

{\bf{Lemma 3}}. Let $X$, $Y$, $Z$ be the feet of the internal bisectors of the angles $\angle A$, $\angle B$, $\angle C$ and let $U$, $V$, $W$ be the points on the sidelines of the triangle $ABC$ defined in Lemma 2. Let the circles $(AXU)$, $(BYV)$, $(CZW)$ meet the triangle circumcircle $\mathcal O$ again at points $A^{\prime}$, $B^{\prime}$, $C^{\prime}$. The lines $AA^{\prime}$, $BB^{\prime}$, $CC^{\prime}$ concur at a point $Q$, which lies on the triangle circumhyperbola $\mathcal H = (IG)^{*}$, the isogonal conjugate of the incenter-centroid line $IG$.

\bigskip

{\it{Proof of Lemma 3}}. By Lemma 2, the points $U$, $V$, $W$ are collinear. By Lemma 1, the radical axes $AA^{\prime}$, $BB^{\prime}$, $CC^{\prime}$ of the triangle circumcircle $\mathcal O$ with the circles $(AXU)$, $(BYV)$, $(CZW)$ concur at a point $Q$. Let $X^{\prime}$, $Y^{\prime}$, $Z^{\prime}$ be the feet of the external bisectors of the angles $\angle A$, $\angle B$, $\angle C$ and $a = \overline{BC}$, $b = \overline{CA}$, $c = \overline{AB}$ the side lengths of the triangle $ABC$. Let $k$ be an arbitrary similarity coefficient of the centrally similar triangles $A_{0}B_{0}C_{0}$, $ABC$, with similarity center at their common incenter $I$. Since $\overline{IA_{0}} = k\ \overline{IA}$ and the external bisectors $A_{0}U$, $AX^{\prime}$ of the angles $\angle A_{0},$ $\angle A$ are parallel, 

$$\frac{\overline{UX^{\prime}}}{\overline{XX^{\prime}}} = \frac{\overline{A_{0}A}}{\overline{XA}} = (1 - k) \frac{\overline{IA}}{\overline{XA}} = \frac{(1 - k)(b + c)}{a + b + c}.$$ 

Using this, the harmonic cross ratio $\frac{\overline{XB}}{\overline{XC}} = -\frac{\overline{X^{\prime}B}}{\overline{X^{\prime}C}} = - \frac{c}{b}$, the ratio $\frac{\overline{XB}}{\overline{X^{\prime}B}} = \frac{c - b}{c + b}$ and Lemma 1, we obtain

$$\frac{\overline{UB}}{\overline{UC}} = \frac{\overline{UX^{\prime}} + \overline{X^{\prime}B}}{\overline{UX^{\prime}} + \overline{X^{\prime}C}} =\frac{c}{b} \cdot \frac{c + a + (2k - 1)b}{a + b + (2k - 1)c}.$$

$$\frac{\overline{Q_{a}B}}{\overline{Q_{a}C}} = \frac{\overline{XB}}{\overline{XC}} \cdot \frac{\overline{UB}}{\overline{UC}} = -\frac{c^{2}}{b^{2}} \cdot \frac{c + a + (2k-1)b}{a + b + (2k-1)c}.$$

With cyclic exchange, this yields the barycentric coordinates of the point $Q$ and its isogonal conjugate $Q^{*}$:

$$Q = \left(\frac{a^{2}}{b + c + (2k-1)a} : \frac{b^{2}}{c + a + (2k-1)b} : \frac{c^{2}}{a + b + (2k-1)c}\right),$$
 
$$Q^{*} = \left(b + c + (2k-1)a : c + a + (2k-1)b : a + b + (2k-1)c\right).$$ 

For $k = 0$, $\frac{1}{2}$, $1$ and $k \rightarrow \infty$, we obtain the four notable points on the line $IG$: the Nagel point $M$, the Spiker point $S = I_{C}$ (the complement of the incenter), the centroid $G$ and the incenter $I$.

Their isogonal conjugates with respect the triangle $ABC$ are the exsimilicenter of the circumcircle and incircle, known as the Kimberling center $X_{56} = M^{*}$ [4], the common point of the Brocard axes of the triangles $IBC$, $ICA$, $IAB$, $ABC$ [3], known as the Kimberling center $X_{58} = S^{*}$, the symmedian point $K = G^{*}$ and the incenter $I = I^{*}$. 

For the points $Q$, $Q^{*}$, corresponding to an arbitrary similarity coefficient $k$ of the triangles $A_{0}B_{0}C_{0}$, $ABC$, the barycentric equation of the line $GQ^{*}$ is $(\vec{G} \times \vec{Q}^{*}) \cdot \vec{R} = 0,$ where $\vec{G},$ $\vec{Q}^{*}$, $\vec{R}$ are 3-dimmensional vectors with components equal to the barycentric coordinates of the points $G$, $Q^{*}$ and an arbitrary point on this line $R = (x : y : z)$, or

$$2(1 - k) \left((b - c)x + (c - a)y + (a - b)z\right) = 0.$$

Since the equation of the line $GQ^{*}$ does not depend on $k$, the point $Q^{*}$ lies on the line $IG$ for any $k$ and the point $Q$ lies on the circumhyperbola $\mathcal H = (IG)^{*}$. $\square$

\bigskip

{\bf{Lemma 4}}. Let the Apollonius circles of the triangle $IBC$ through the vertices $B$, $C$ cut the circumcircle $\mathcal O$ of the reference triangle $ABC$ again at points $B^{\prime}$, $C^{\prime}$. Let the lines $BB^{\prime}$, $CC^{\prime}$ intersect at a point $Q_{BC}$. Let $Q_{CA}$, $Q_{AB}$ be analogously defined points for the triangles $ICA$, $IAB$. The points $Q_{BC}$, $Q_{CA}$, $Q_{AB}$ lie on the circumhyperbola $\mathcal H = (IG)^{*}$ and their barycentric coordinates are obtained from Lemma 3 with the similarity coefficients $k = -\cos A$, $-\cos B$, $-\cos C$, respectively.

\bigskip

{\it{Proof of Lemma 4}}. Let $X$, $Y$, $Z$ be the feet of the internal bisectors and $X^{\prime}$, $Y^{\prime}$, $Z^{\prime}$ the feet of the external bisectors of the angles $\angle A$, $\angle B$, $\angle C$ of the triangle $ABC$. Let its Apollonius circle $(AXX^{\prime})$ through the vertex $A$ cut the lines $CA$, $AB$ again at points $B_{a}$, $C_{a}$. Since $AX'$ bisects the angle $\angle C_{a}AB_{a}$, $X'$ is the midpoint of the arc $B_{a}C_{a}$ of the circumcircle $(AXX^{\prime})$ of the triangle $AB_{a}C_{a}$, hence, $C_{a}$ is the reflection of $B_{a}$ in the line $BC$.

Let $\mathcal S_{b}$, $\mathcal S_{c}$ be the Apollonius circles of the triangle $IBC$ through the vertices $B$, $C$. Since $CA$ is the reflection of $CB$ in $CI$, the circle $\mathcal S_{b}$ cuts the line $BI$ again at the point $Y$ and similarly, the circle $\mathcal S_{c}$ cuts the line $CI$ again at the point $Z$. The circumcircle $\mathcal M$ of the triangle $IBC$ is centered at the midpoint $M^{\prime}$ of the arc $BC$ of the circumcircle $\mathcal O$ of the reference triangle $ABC$ opposite to the vertex $A$. The centers $S_{b}$, $S_{c}$ of the Apollonius circles $\mathcal S_{b}$, $\mathcal S_{c}$ are the intersections of the lines $CI$, $BI$ with the tangents to the circle $\mathcal M$ at the points $B$, $C$. Let the circles $\mathcal S_{b}$, $\mathcal S_{c}$ cut the lines $CA$, $AB$ again at points $V$, $W$. The triangles $BCV$, $CBW$ are both isosceles and consequently,

$$\angle IBV = \angle ICW = \frac{\pi}{2} - \left(\frac{\angle B}{2} + \frac{\angle C}{2}\right) = \frac{\angle A}{2}.$$

Let $B_{0}$, $C_{0}$ be the feet of the altitudes through $V$, $W$ of the isosceles triangles $BIV$, $CIW$ with the base angles equal to $\frac{\angle A}{2}$. Then,
 
$$\frac{\overline{IB_{0}}}{\overline{IB}} = \frac{\overline{IC_{0}}}{\overline{IC}} = -\cos A,$$

which means that the lines $B_{0}C_{0}$, $BC$ are parallel. Let $A_{0}$ be a point on the internal bisector $AX$ of the angle $\angle A$, such that the triangles $A_{0}B_{0}C_{0}$, $ABC$ are centrally similar with similarity center $I$; their similarity coefficient is equal to $k = -\cos A$. The lines $B_{0}V$, $C_{0}W$, perpendicular to the lines $BI$, $CI$, are the external bisectors of the angles $\angle B_{0}$, $\angle C_{0}$, Let the external bisector of the angle $\angle A_{0}$ intersect the line $BC$ at a point $U$ and let the line $AQ_{BC}$ cut the circumcircle $\mathcal O$ of the triangle $ABC$ again at a point $A^{\prime}$. By Lemma 2, the points $U$, $V$, $W$ are collinear. By definition, the points $B$, $Y$, $B^{\prime}$, $V$ are concyclic and the points $C$, $Z$, $C^{\prime}$, $W$ are concyclic. The points $A$, $X$, $A^{\prime}$, $U$ are then also concyclic, because the opposite would contradict Lemma 1. By Lemma 3, the point $Q_{BC}$ lies on the circumhyperbola $\mathcal H = (IG)^{*}$ and the barycentric coordinates of the point $Q_{BC}$ are obtained from Lemma 3 with the similarity coefficient $k = -\cos A$. $\square$
  
\bigskip

For the incenter $I = I^{*}$, the circles $(AXA^{\prime})$, $(BYB^{\prime})$, $(CZC^{\prime})$ degenerate to a pencil of concurrent lines, the internal bisectors $AX$, $BY$, $CZ$ of the angles $\angle A$, $\angle B$, $\angle C$. For the symmedian point $K = G^{*}$, the circles $(AXA^{\prime})$, $(BYB^{\prime})$, $(CZC^{\prime})$ are the Apollonius circles of the triangle $ABC$, which are coaxal and the Kimberling center $X_{58}$ lies on their common radical axis [3]. For the exsimilicenter of the circumcircle and incircle $X_{56} = M^{*}$, the circles $(AXA^{\prime})$, $(BYB^{\prime})$, $(CZC^{\prime})$ are also coaxal and the Kimberling center $X_{58}$ also lies on their common radical axis [6]. Based on the preliminary results, we now conjecture our main theorem:

\bigskip
\begin{center}
\includegraphics[scale=1.1]{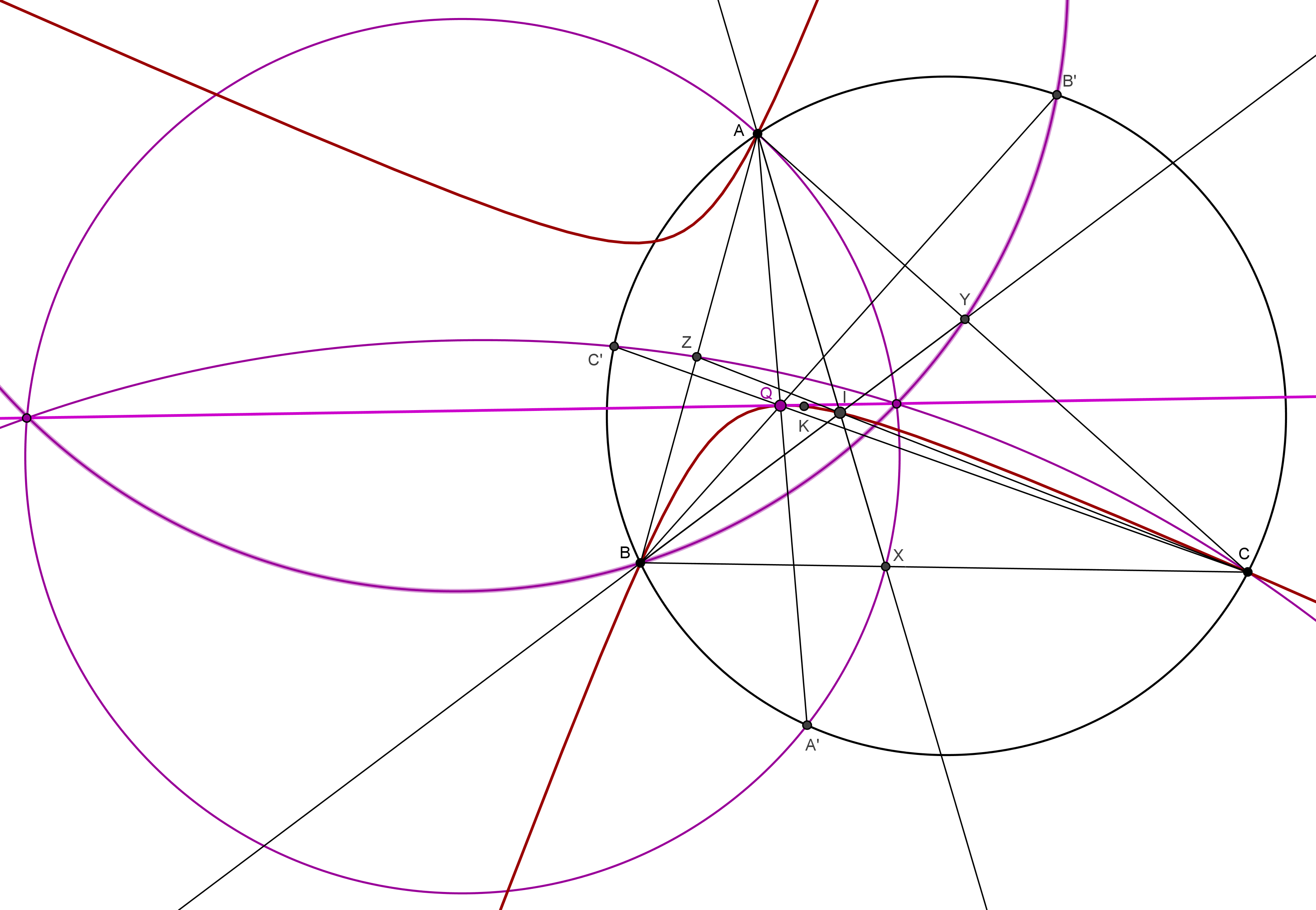}
\end{center}
\bigskip

{\bf{Theorem 5}}. Let $Q$ be an arbitrary point in the plane of a non-equilateral triangle $ABC$, different from the triangle vertices. Let $XYZ$ be the cevian triangle of the incenter $I$ and $A^{\prime}B^{\prime}C^{\prime}$ the circumcevian triangle of the point $Q$. The circles $(AXA^{\prime})$, $(BYB^{\prime})$, $(CZC^{\prime})$ are coaxal and their common radical axis passes through the Kimberling center $X_{58}$, if and only if the point $Q$ lies on the triangle circumhyperbola $\mathcal H = (IG)^{*}$, the isogonal conjugate of the incenter-centroid line $IG$.

\bigskip

{\it{Remark}}. When the point $Q$ in Theorem 5 is identical with the point $Q_{BC}$ from Lemma 4, the two Apollonius circles $(BYB^{\prime})$, $(CZC^{\prime})$ of the triangle $IBC$ through its vertices $B$, $C$, supplemented by the circle $(AXA^{\prime})$, become a special case of Theorem 5. By Theorem 5, the Brocard axis of the triangle $IBC$, identical with the radical of its Apollonius circles $(BYB^{\prime})$, $(CZC^{\prime})$, passes through the Kimberling center $X_{58}$. Similarly, the Brocard axes of the triangles $ICA$, $IAB$ pass through $X_{58}$.

\bigskip

We will postpone the proof of Theorem 5 and at first, analyze the four special cases $X_{56} = M^{*}$, $X_{58} = S^{*}$, $K = G^{*}$ and $I = I^{*}$.

\bigskip

{\bf{2. Special cases}}

\bigskip

Inversion in any circle carries a pencil of circles or lines into a pencil of circles or lines. Hereby, we will use the inversion in the triangle incircle $\mathcal I$ and we will show that the inversion images of the circles $(AXA^{\prime})$, $(BYB^{\prime})$, $(CZC^{\prime})$ for $X_{56} = M^{*}$, $X_{58} = S^{*}$, $K = G^{*}$ and $I = I^{*}$ are coaxal.

\bigskip
 
Consider the inversion $\Psi$ with center $I$ and power $r^2$, where $r$ is the inradius length. Let $D$, $E$, $F$ be the tangency points of the incircle $\mathcal I$ with the triangle sides $BC$, $CA$, $AB$. The inversion $\Psi$ carries the triangle sidelines $BC$, $CA$, $AB$ into the circles $\Gamma_{a}$, $\Gamma_{b}$, $\Gamma_{c}$ tangent to $BC$, $CA$, $ AB$ at the points $D$, $E$, $F$ and passing through the inversion center $I$. The tangency points $D$, $E$, $F$, lying on the inversion circle, are invariant under $\Psi$. Having equal diameters $ID = IE = IF = r$, these three circles are congruent. They pairwise intersect at the midpoints $A_{1}$, $B_{1}$, $C_{1}$ of the sides $EF$, $FD$, $DE$ of the contact triangle $DEF$, the images of the vertices $A$, $B$, $C$ of the reference triangle $ABC$. The image $\mathcal O_{1}$ of the circumcircle $\mathcal O$ of the triangle $ABC$ is the circumcircle of the medial triangle $A_{1}B_{1}C_{1}$ of the contact triangle $DEF$. Thus, $\mathcal O_{1}$ is the nine-point circle of the contact triangle with radius $\frac{r}{2}$, congruent to the circles $\Gamma_{a}$, $\Gamma_{b}$, $\Gamma_{c}$. The internal angle bisectors $AX$, $BY$, $CZ$, passing through the inversion center $I$, are carried into themselves. The images $X_{1}$, $Y_{1}$, $Z_{1}$ of their feet $X$, $Y$, $Z$ are therefore the intersections of the lines $AX$, $BY$, $CZ$ with the circles $\Gamma_{a}$, $\Gamma_{b}$, $\Gamma_{c}$ other than $I$. Since $I$ is the orthocenter of the triangle $A_{1}B_{1}C_{1}$ and the circles $\Gamma_{a}$, $\Gamma_{b}$, $\Gamma_{c}$ are the sideline reflections of its circumcircle $\mathcal O_{1}$, the points $X_{1}$, $Y_{1}$, $Z_{1}$ are the sideline reflections of its vertices $A_{1}$, $B_{1}$, $C_{1}$. The images $(A_{1}X_{1}A_{1}^{\prime})$, $(B_{1}Y_{1}B_{1}^{\prime})$, $(C_{1}Z_{1}C_{1}^{\prime})$ of the circles $(AXA^{\prime})$, $(BYB^{\prime})$, $(CZC^{\prime})$ under $\Psi$ pass through the vertices $A_{1}$, $B_{1}$, $C_{1}$ of the inverted triangle $A_{1}B_{1}C_{1}$ and through their sideline reflections $X_{1}$, $Y_{1}$, $Z_{1}$, which means that they are centered on the respective sidelines of this triangle.

\bigskip

{\bf{2. 1. The exsimilicenter of the circumcircle and incircle}}

\bigskip

When the point $Q$ in Theorem 5 is identical with the exsimilicenter of the circumcircle and incircle, the Kimberling center $X_{56} = M^{*}$ and isogonal conjugate of the Nagel point $M$, the circles $(AXA^{\prime})$, $(BYB^{\prime})$, $(CZC^{\prime})$ intersect the circumcircle $\mathcal O$ at its tangency points $A^{\prime}$, $B^{\prime}$, $C^{\prime}$ with the triangle mixtilinear incircles in the angles $\angle A$, $\angle B$, $\angle C$ [7]. Indeed, the tangency point $A^{\prime}$ is the external similarity center of the circumcircle $\mathcal O$ and mixtilinear incircle $\mathcal K_{a}$ in the angle $\angle A$, while the triangle vertex $A$ is the external similarity center of the incircle $\mathcal I$ and mixtilinear incircle $\mathcal K_{a}$. By the Monge-d'Alembert theorem, the external similarity center $X_{56}$ of the circles $\mathcal O$, $\mathcal I$ lies on the line $AA^{\prime}$ and similarly, $X_{56}$ lies on the lines $BB^{\prime}$, $CC^{\prime}$ [5]. 
 
Consider the inversion $\Psi$ in the incircle $\mathcal I$ defined previously. Since the congruent circles $\mathcal O_{1},$ $\Gamma_{b},$ $\Gamma_{c}$ with radii $\frac{r}{2}$ meet at the point $A_{1}$, a circle $\mathcal L_{a}$ with radius $r$ (congruent to the incircle $\mathcal I$) and centered at $A_{1}$ and is tangent to all three at antipodal points with respect to $A_{1}$. Four different circles are tangent to the triangle sidelines $AB$ $AC$ and to the circumcircle $\mathcal O$: the mixtilinear incircle $\mathcal K_{a}$ and mixtilinear excircle $\mathcal K_{a}^{\prime}$ in the angle $\angle A$, both centered on the ray $AX$ of the internal bisector of the angle $\angle A$, and two other circles centered on the opposite rays of the external bisector of this angle. Since the inversion center $I$ is the similarity center of a circle and its inversion image, only the images of the mixtilinear incircle $\mathcal K_{a}$ and mixtilinear excircle $\mathcal K_{a}^{\prime}$ are centered on the internal angle bisector $AX$. The mixtilinear excircle lies outside the circumcircle $\mathcal O$ and outside the inversion circle $\mathcal I$. Therefore, it has both intersections with the line $AX$ on the ray $IX$. Since the inversion $\Psi$ has positive power $r^2$, its image also has both intersections with the line $AX$ on the ray $IX$. Consequently, it is centered on the ray $IX$ and it cannot be identical with the circle $\mathcal L_{a}$ centered on the opposite ray $IA$. As a result, the circle $\mathcal L_{a}$ is the inversion image of the mixtilinear incircle $\mathcal K_{a}$ in the angle $\angle A$ and similarly, the circles $\mathcal L_{b}$, $\mathcal L_{c}$ with radii $r$ and centered at the points $B_{1}$, $C_{1}$ are the images of the mixtilinear incircles $\mathcal K_{b}$, $\mathcal K_{c}$ in the angles $\angle B$, $\angle C$. Thus, the images $A_{1}^{\prime}$, $B_{1}^{\prime}$, $C_{1}^{\prime}$ of the tangency points $A^{\prime}$, $B^{\prime}$, $C^{\prime}$ of the mixtilinear incircles $\mathcal K_{a}$, $\mathcal K_{b}$, $\mathcal K_{c}$ with the circumcircle $\mathcal O$ are the antipodal points of the circumcircle $\mathcal O_{1}$ of the triangle $A_{1}B_{1}C_{1}$ with respect to its vertices $A_{1}$, $B_{1}$, $C_{1}$.

The images $(A_{1}X_{1}A_{1}^{\prime})$, $(B_{1}Y_{1}B_{1}^{\prime})$, $(C_{1}Z_{1}C_{1}^{\prime})$ of the circles $(AXA^{\prime})$, $(BYB^{\prime})$, $(CZC^{\prime})$ under $\Psi$ are centered on the respective sidelines of the inverted triangle $A_{1}B_{1}C_{1}$ and they pass through the endpoints of the diameters $A_{1}A_{1}^{\prime}$, $B_{1}B_{1}^{\prime}$, $C_{1}C_{1}^{\prime}$ of its circumcircle $\mathcal O_{1}$. Thus, their centers $O_{a}$, $O_{b}$, $O_{c}$ are the intersections of the perpendicular bisectors of the circumcircle diameters $A_{1}A_{1}^{\prime}$, $B_{1}B_{1}^{\prime}$, $C_{1}C_{1}^{\prime}$ through their common midpoint $O_{1}$ with the respective triangle sidelines. Consequently, their centers $O_{a}$, $O_{b}$, $O_{c}$ are collinear and their pairwise radical axes parallel, the line $O_{a}O_{b}O_{c}$ being the orthotransversal of the triangle $A_{1}B_{1}C_{1}$ with respect to its circumcenter $O_{1}$ [2], which lies on the Euler line $O_{1}I$. Powers of $O_{1}$ to the circles  $(A_{1}X_{1}A_{1}^{\prime})$, $(B_{1}Y_{1}B_{1}^{\prime})$, $(C_{1}Z_{1}C_{1}^{\prime})$ are equal to the power of $O_{1}$ to the circumcircle $\mathcal O_{1}$. Existence of the point $O_{1}$ with equal powers to all three circles then implies a common radical axis, {\it{i.e.}}, the circles $(A_{1}X_{1}A_{1}^{\prime})$, $(B_{1}Y_{1}B_{1}^{\prime})$, $(C_{1}Z_{1}C_{1}^{\prime})$ are coaxal and the circumcenter $O_{1}$ lies on their common radical axis.  Thus, the coaxality claim of Theorem 5 holds for $X_{56} = M^{*}$. $\square$

\bigskip

Additional results, easily obtained from the inverted figure, are collected in the following lemma:

\bigskip
\begin{center}
\includegraphics[scale=0.35]{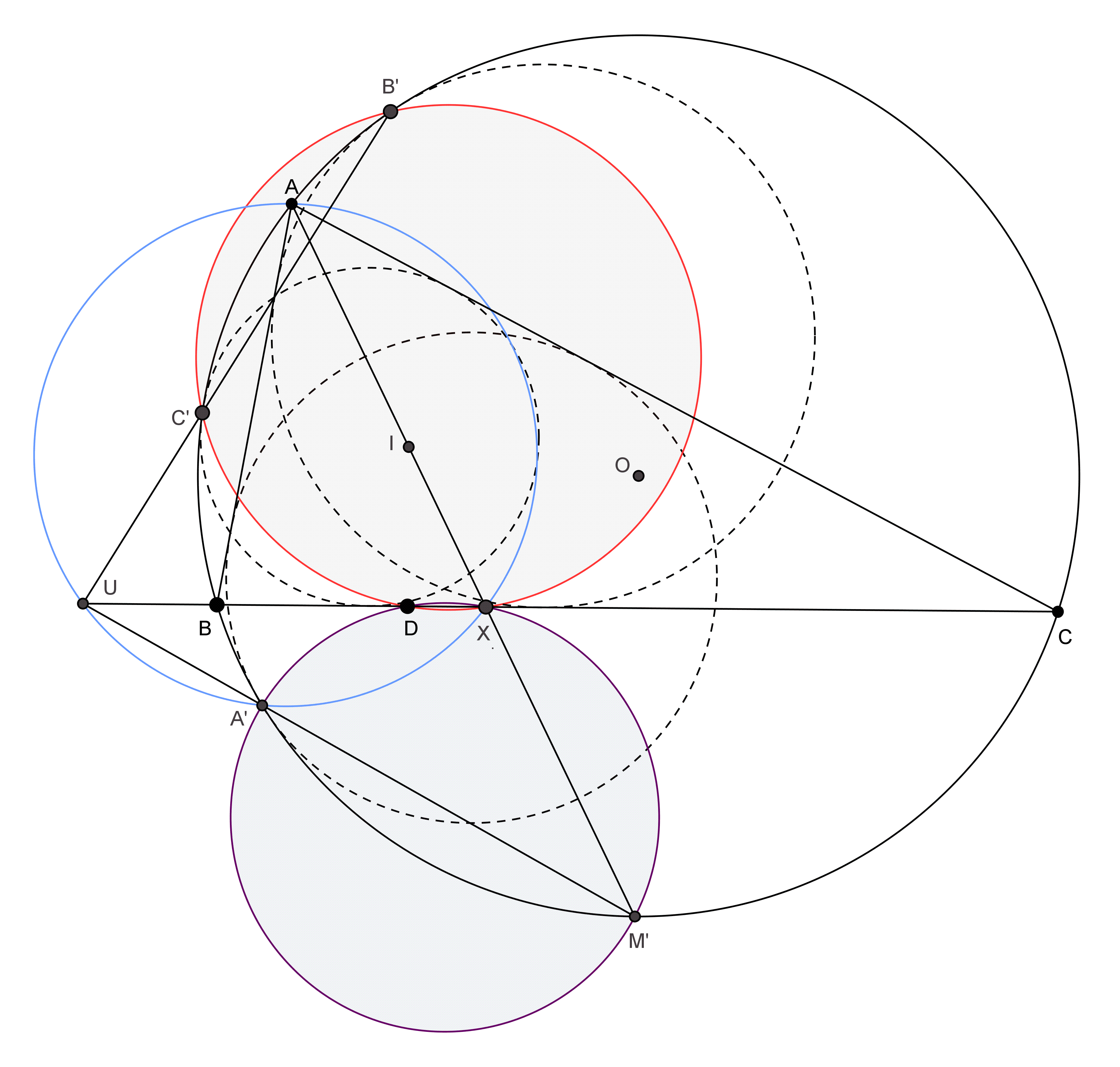}
\end{center}
\bigskip

{\bf{Lemma 6}}. Let $M^{\prime}$, $N^{\prime}$, $P^{\prime}$ be the midpoints of the circumcircle arcs $BC$, $CA$, $AB$ opposite to the triangle vertices $A$, $B$, $C$. The quadrilaterals $A^{\prime}M^{\prime}XD$, $B^{\prime}N^{\prime}YE$, $C^{\prime}P^{\prime}ZF$ are cyclic and the quadrilaterals $B^{\prime}C^{\prime}XD$, $C^{\prime}A^{\prime}YE$, $A^{\prime}B^{\prime}ZF$ are also cyclic. The radical centers $U$, $V$, $W$ of the circle pairs $(A^{\prime}M^{\prime}XD)$, $(B^{\prime}C^{\prime}XD)$; $(B^{\prime}N^{\prime}YE)$, $(C^{\prime}A^{\prime}YE)$; $(C^{\prime}P^{\prime}ZF)$, $(A^{\prime}B^{\prime}ZF)$ with the circumcircle $\mathcal O$ are collinear and they lie on the circles $(AXA^{\prime})$, $(BYB^{\prime})$, $(CZC^{\prime})$, respectively.

\bigskip

{\it{Proof of Lemma 6}}. The images $M_{1}^{\prime}$, $N_{1}^{\prime}$, $P_{1}^{\prime}$ of the points $M^{\prime}$, $N^{\prime}$, $P^{\prime}$ under the inversion $\Psi$ are the intersections of the circle $\mathcal O_{1}$ with the internal angle bisectors $AX$, $BY$, $CZ$ other than the points $A_{1}$, $B_{1}$, $C_{1}$. Since $A_{1}A_{1}^{\prime}$ is a diameter of $\mathcal O_{1}$, the angle $\angle IM_{1}^{\prime}A_{1}^{\prime} = \angle A_{1}M_{1}^{\prime}A_{1}^{\prime}$ is right. By the basic properties of the inversion $\Psi$, the triangles $IA^{\prime}M^{\prime}$, $IM_{1}^{\prime}A_{1}^{\prime}$ are oppositely similar and the angle $\angle M^{\prime}A^{\prime}I = \angle IM_{1}^{\prime}A_{1}^{\prime}$ is also right. Similarly, the angles $\angle N^{\prime}B^{\prime}I$, $\angle P^{\prime}C^{\prime}I$ are right. Let the lines $M^{\prime}A^{\prime}$, $N^{\prime}B^{\prime}$, $P^{\prime}C^{\prime}$ cut the triangle sidelines $BC$, $CA$, $AB$ at points $U$, $V$, $W$. The line $BC$ is the image of the circumcircle $\mathcal O$ in the inversion in the circumcircle $\mathcal M$ of the triangle $IBC$, centered at the point $M^{\prime}$, and the point $U$ is the image of the point $A^{\prime}$ in this inversion. The line $IA^{\prime}$, being perpendicular to the line $M^{\prime}U$ at $A^{\prime}$, is the polar of the point $U$ with respect to the circle $\mathcal M$. Consequently, the line $IU$ is the tangent of the circle $\mathcal M$ at the point $I$, perpendicular to the internal bisector $AM^{\prime}$ of the angle $\angle A$. Similarly, the lines $IV$, $IW$ are perpendicular to the internal bisectors $BN^{\prime}$, $CP^{\prime}$ of the angles $\angle B$, $\angle C$. This means the points $U$, $V$, $W$ are collinear, the line $UVW$ being the orthotransversal with respect to the incenter $I$ [2].

The image $U_{1}$ of the point $U$ under $\Psi$ is the intersection of the line $IU$ with the circle $\Gamma_{a}$ other than $I$, hence, the quadrilateral $IX_{1}DU_{1}$ is a rectangle inscribed in the circle $\Gamma_{a}$. The right angle triangles $IX_{1}D$, $A_{1}M_{1}^{\prime}A_{1}^{\prime}$ with the common sideline  $IX_{1} = A_{1}M_{1}^{\prime}$ inscribed in the congruent circles $\Gamma_{a}$, $\mathcal O_{1}$ are congruent. It follows that the quadrilateral $A_{1}^{\prime}M_{1}^{\prime}X_{1}D$ is also a rectangle and the quadrilateral $A_{1}X_{1}A_{1}^{\prime}U_{1}$ is an isosceles trapezoid, both cyclic. Consequently, the quadrilaterals $A^{\prime}M^{\prime}XD$ and $AXA^{\prime}U$ are also cyclic. Since the segments $B_{1}B_{1}^{\prime}$, $C_{1}C_{1}^{\prime}$ are diameters of the circle $\mathcal O_{1}$, the quadrilateral $B_{1}C_{1}B_{1}^{\prime}C_{1}^{\prime}$ is a rectangle and the lines $B_{1}C_{1}$, $C_{1}^{\prime}B_{1}^{\prime}$ are parallel. The midline $B_{1}C_{1}$ of the contact triangle $DEF$ is parallel to its sideline $EF$. The lines $X_{1}D$, $EF$ are also parallel, being both perpendicular to the angle bisector $AM^{\prime}$. Combining, the three lines $B_{1}C_{1}$, $C_{1}^{\prime}B_{1}^{\prime}$, $X_{1}D$ are parallel. But, $C_{1}^{\prime}B_{1}^{\prime}$, $X_{1}D$ are chords of the circles $\mathcal O_{1}$, $\Gamma_{a}$, parallel to their radical axis $B_{1}C_{1}$. It follows that the quadrilateral $B_{1}^{\prime}C_{1}^{\prime}X_{1}D$ is an isosceles trapezoid, which is cyclic, and consequently, the quadrilateral $B^{\prime}C^{\prime}XD$ is also cyclic. Moreover, the pairwise radical axes $M^{\prime}A^{\prime}$, $C^{\prime}B^{\prime}$, $BC$ of the circles $\mathcal O$, $(A^{\prime}M^{\prime}XD)$, $(B^{\prime}C^{\prime}XD)$ meet at their radical center $U$. $\square$

\bigskip

{\it{Remark}}. The orthotransversal $UVW$ of the triangle $ABC$ with respect to the incenter $I$ is the tripolar with respect to its orthocorrespondent $I^{\perp}$ [2]. Using Lemma 3 with the similarity coefficient $k = 0$, we obtain the barycentric coordinates of $I^{\perp}$:

$$I^{\perp} = \left(\frac{a}{b + c - a} : \frac{b}{c + a - b} : \frac{c}{a + b - c}\right).$$

Thus, $I^{\perp}$ is identical with the Kimberling center $X_{57}$. Neither $I^{\perp}$, nor the identity $I^{\perp} = X_{57}$ are included in the current edition of [4].

\bigskip

{\bf{2. 2. The common point of the Brocard axes of triangles $IBC$, $ICA$, $IAB$, and $ABC$}}

\bigskip

Assume that the point $Q$ in Theorem 5 is identical with the common point of the Brocard axes of the triangles $IBC$, $ICA$, $IAB$, $ABC$ [3], the Kimberling center $X_{58} = S^{*}$ and isogonal conjugate of the Spiker point $S$. By Lemma 3 with the similarity coefficient $k = \frac{1}{2}$, the circles $(AXA^{\prime})$, $(BYB^{\prime})$, $(CZC^{\prime})$ pass through the collinear intersections $U$, $V$, $W$ of the perpendicular bisectors of the angle bisector segments $IA$, $IB$, $IC$ with the triangle sidelines $BC$, $CA$, $AB$. 

Let $A_{0}$, $B_{0}$, $C_{0}$ be the midpoints of the segments $IA$, $IB$, $IC$. The inversion $\Psi$ in the incircle $\mathcal I$ carries the perpendicular bisectors $A_{0}U,$ $B_{0}V$, $C_{0}W$ of these segments into the circles $\mathcal U_{1}$, $\mathcal V_{1}$, $\mathcal W_{1}$ centered on the lines $AI$, $BI$ $CI$ and passing through the inversion center $I$ (Figure 4). Let $r_{u}$ $r_{v}$, $r_{w}$ be radii of these circles, respectively. By the basic properties of the inversion $\Psi$, we have

$$r^{2} = 2r_{u} \cdot IA_{0} = r_{u} \cdot IA.$$

On the other hand, $EA_{1}$ is the altitude through $E$ of the right angle triangle $IAE$, hence,

$$\frac{IA}{IE} = \frac{IE}{IA_{1}},\ IA_{1} \cdot IA = IE^{2} = r^2.$$

As a result, $r_{u} = IA_{1}$ and similarly, $r_{v} = IB_{1}$, $r_{w} = IC_{1}$, which means that the circles $\mathcal U_{1}$, $\mathcal V_{1}$, $\mathcal W_{1}$ are centered at the vertices $A_{1}$, $B_{1}$, $C_{1}$ of the inverted triangle $A_{1}B_{1}C_{1}$. The images $U_{1}$, $V_{1}$, $W_{1}$ of the points $U$, $V$, $W$ under $\Psi$ are the intersections of the circle pairs $\mathcal U_{1}$, $\Gamma_{a}$; $\mathcal V_{1}$, $\Gamma_{b}$; $\mathcal W_{1}$, $\Gamma_{c}$ other than $I$.  

Let $D_{1}$, $E_{1}$, $F_{1}$ be the centers of the circles $\Gamma_{a}$, $\Gamma_{b}$, $\Gamma_{c}$ with diameters $ID = IE = IF = r$. Since $\Gamma_{a}$ is the reflection of the circumcircle $\mathcal O_{1}$ of the triangle $A_{1}B_{1}C_{1}$ in its sideline $B_{1}C_{1}$, the quadrilateral $A_{1}O_{1}D_{1}I$ is a parallelogram and $D_{1}O_{1} = IA_{1}$. Let the line $D_{1}O_{1}$ parallel to $IA_{1}$ meet the arc $B_{1}C_{1}$ of the circumcircle $\mathcal O_{1}$ opposite to $D_{1}$ at a point $P_{1}$. Again, the quadrilateral $A_{1}P_{1}U_{1}I$ is a parallelogram and in addition, $U_{1}P_{1} = U_{1}A_{1} = IA_{1}$. Because the triangle $A_{1}U_{1}P_{1}$ is isosceles, the line $O_{1}U_{1}$ is the perpendicular bisector of the circumcircle chord $A_{1}P_{1}$. On the other hand, the center line $D_{1}A_{1}$ of the circles $\Gamma_{a}$, $\mathcal U_{1}$ is the perpendicular bisector of the radical axis segment $IU_{1}$. Since the perpendicular bisectors $O_{1}U_{1}$, $D_{1}A_{1}$ of the parallel segments $A_{1}P_{1}$, $IU_{1}$ are parallel and since $D_{1}O_{1} = U_{1}A_{1}$ the quadrilateral $A_{1}D_{1}O_{1}U_{1}$ is an isosceles trapezoid, and similarly, the quadrilaterals $B_{1}E_{1}O_{1}V_{1}$, $C_{1}F_{1}O_{1}W_{1}$ are also isosceles trapezoids, all cyclic. The circumcircles $\mathcal O_{a}$, $\mathcal O_{b}$, $\mathcal O_{c}$ of these isosceles trapezoids are centered on the perpendicular bisectors of their sides $D_{1}O_{1}$, $E_{1}O_{1}$, $F_{1}O_{1}$, identical with the respective sidelines of the triangle $A_{1}B_{1}C_{1}$. It follows that the images $X_{1}$, $Y_{1}$, $Z_{1}$ of the points $X$, $Y$, $Z$ under $\Psi$, identical with the reflections of the triangles vertices $A_{1}$, $B_{1}$, $C_{1}$ in its sidelines, also lie on the circles $\mathcal O_{a}$, $\mathcal O_{b}$, $\mathcal O_{c}$. Hereby, these three circles are the images of the circles $(AXA^{\prime}U)$  $(BYB^{\prime}V)$, $(CZC^{\prime}W)$ under $\Psi$. The images $A_{1}^{\prime}$, $B_{1}^{\prime}$, $C_{1}^{\prime}$ of the points $A^{\prime}$, $B^{\prime}$ $C^{\prime}$ are then their intersections with the circumcircle $\mathcal O_{1}$ other than $A_{1}$, $B_{1}$, $C_{1}$, identical with the reflections of the points $U_{1}$, $V_{1}$, $W_{1}$ in the respective triangle sidelines.

The nine-point circle center $N_{1}$ of the triangle $A_{1}B_{1}C_{1}$ is the common midpoint of the chords $A_{1}D_{1}$, $B_{1}E_{1}$, $C_{1}F_{1}$ of the circles $\mathcal O_{a}$, $\mathcal O_{b}$, $\mathcal O_{c}$. Thus, their centers $O_{a}$, $O_{b}$, $O_{c}$ are the intersections of the perpendicular bisectors of the chords $A_{1}D_{1}$, $B_{1}E_{1}$, $C_{1}F_{1}$ through their common midpoint $N_{1}$ with the respective triangle sidelines. Consequently, the centers $O_{a}$, $O_{b}$, $O_{c}$ are collinear, the line $O_{a}O_{b}O_{c}$ being the orthotransversal of the triangle $A_{1}B_{1}C_{1}$ with respect to its nine-point circle center $N_{1}$ [2], which lies on the Euler line $O_{1}I$. In addition, the circles $\mathcal O_{a}$, $\mathcal O_{b}$, $\mathcal O_{c}$ concur at the circumcenter $O_{1}$, hence, powers of $O_{1}$ to the circles $\mathcal O_{a}$, $\mathcal O_{b}$, $\mathcal O_{c}$ are zero. Existence of the point $O_{1}$ with equal powers to all three circles then implies a common radical axis, {\it{i.e.}}, the circles $\mathcal O_{a}$, $\mathcal O_{b}$, $\mathcal O_{c}$ are coaxal with and the circumcenter $O_{1}$ lies on their common radical axis. Thus, the coaxality claim of Theorem 5 holds for $X_{58} = S^{*}$. $\square$
 
\bigskip

{\bf{2. 3. The symmedian point}}

\bigskip

When the point $Q$ in Theorem 5 is identical with the symmedian point $K = G^{*}$, the isogonal conjugate of the centroid $G$, the circles $(AXA^{\prime})$, $(BYB^{\prime})$, $(CZC^{\prime})$ are the Apollonius circles of the triangle $ABC$ [1]. Indeed, the symmedian $AA^{\prime}$ passes through the intersection $T_{a}$ of the circumcircle tangents at the triangle vertices $B$, $C$. Inversion in the circumcircle $\mathcal O$ carries the line $AA^{\prime}T_{a}$ into the circle $(AOM_{a}A')$ passing through the inversion center $O$ and through the midpoint $M_{a}$ of the triangle side $BC$, the image of $T_{a}$. The circle $(AXA^{\prime})$ is centered at the intersection $S_{a}$ of the circumcircle tangents at the points $A$, $A^{\prime}$ and the segment $OS_{a}$ is a diameter of the circle $(AOM_{a}A^{\prime})$. Because of the right angle $\angle OM_{a}S_{a}$, the center $S_{a}$ of the circle $(AXA^{\prime})$ is on the sideline $BC$ and consequently, it is one of the three Apollonius circles of the triangle $ABC$. The Apollonius circles $(AXA^{\prime})$, $(BYB^{\prime})$, $(CZC^{\prime})$ meet at the two isodynamic points of the triangle, hence, they are coaxal. Their centers are the intersections of the circumcircle tangents at the vertices $A$, $B$, $C$ with the sidelines $BC$, $CA$, $AB$ and conversely, $OA$, $OB$, $OC$ are their tangents from the circumcenter $O$. They are therefore perpendicular to the circumcircle $\mathcal O$. Since powers of $O$ to the Apollonius circles are equal, the circumcenter $O$ lies on their common radical axis. Since powers of $K$ to the Apollonius circles $(AXA^{\prime})$, $(BYB^{\prime})$, $(CZC^{\prime})$ are equal to the power of $K$ to the circumcircle $\mathcal O$, the symmedian point $K$ also lies on their common radical axis, the Brocard axis of the triangle $ABC$.

The images $(A_{1}X_{1}A_{1}^{\prime})$, $(B_{1}Y_{1}B_{1}^{\prime})$, $(C_{1}Z_{1}C_{1}^{\prime})$ of the Apollonius circles $(AXA^{\prime})$, $(BYB^{\prime})$, $(CZC^{\prime})$ of the triangle $ABC$ under the inversion $\Psi$ in the incircle $\mathcal I$ are centered on the respective sidelines of the inverted triangle $A_{1}B_{1}C_{1}$ and they pass through its vertices $A_{1}$, $B_{1}$, $C_{1}$. Since the Apollonius circles of the triangle $ABC$ are perpendicular its circumcircle $\mathcal O$, their images are also perpendicular to the circumcircle image $\mathcal O_{1}$. As a result, the images of the Apollonius circles of the triangle $ABC$ are the Apollonius circles of the inverted triangle $A_{1}B_{1}C_{1}$. They are coaxal and the triangle circumcenter $O_{1}$ lies on their common radical axis, the Brocard axis of the triangle $A_{1}B_{1}C_{1}$.  Thus, the coaxality claim of Theorem 5 holds for $K = G^{*}$. $\square$
 
\bigskip

Let $O_{a}$, $O_{b}$, $O_{c}$ be the centers of the Apollonius circles of the triangle $A_{1}B_{1}C_{1}$. Let $\mathcal P_{a}$, $\mathcal P_{b}$, $\mathcal P_{c}$ be circles with diameters $A_{1}O_{a}$, $B_{1}O_{b}$, $C_{1}O_{c}$, tangent to the Apollonius circles at the triangle vertices $A_{1}$, $B_{1}$, $C_{1}$. The circumcircle radii $O_{1}A_{1}$, $O_{1}B_{1}$, $O_{1}C_{1}$ are the common tangents of the circles $\mathcal P_{a}$, $\mathcal P_{b}$, $\mathcal P_{c}$ and the Apollonius circles, hence, powers of the circumcenter $O_{1}$ to the circles $\mathcal P_{a}$, $\mathcal P_{b}$, $\mathcal P_{c}$ are equal,

$$\overline{O_{1}A_{1}}^2 = \overline{O_{1}B_{1}}^2 = \overline{O_{1}C_{1}}^2.$$

Since the circles $\mathcal P_{a}$, $\mathcal P_{b}$, $\mathcal P_{c}$ are centered on the respective midlines of the triangle $A_{1}B_{1}C_{1}$ and pass through its vertices $A_{1}$, $B_{1}$, $C_{1}$, they also pass through its altitude feet $J_{a}$, $J_{b}$, $J_{c}$. Let the altitudes $A_{1}J_{a}$, $B_{1}J_{b}$, $C_{1}J_{c}$ meet the circumcircle $\mathcal O_{1}$ again at points $H_{a}$, $H_{b}$, $H_{c}$, reflections of the orthocenter $I$ in the triangle sidelines $B_{1}C_{1}$, $C_{1}A_{1}$, $A_{1}B_{1}$. Powers of the orthocenter $I$ to the circles $\mathcal P_{a}$, $\mathcal P_{b}$, $\mathcal P_{c}$ are then

$$\overline{IJ_{a}} \cdot \overline{IA_{1}} = \frac{1}{2} \overline{IH_{a}} \cdot \overline{IA_{1}},\ \overline{IJ_{b}} \cdot \overline{IB_{1}} = \frac{1}{2} \overline{IH_{b}} \cdot \overline{IB_{1}},\ \overline{IJ_{c}} \cdot \overline{IC_{1}} = \frac{1}{2} \overline{IH_{c}} \cdot \overline{IC_{1}},$$

respectively. Thus, powers of the orthocenter $I$ to all three circles are the same, equal to half the power of $I$ to the circumcircle $\mathcal O_{1}$. It follows that the circles $\mathcal P_{a}$, $\mathcal P_{b}$, $\mathcal P_{c}$ are coaxal and the Euler line $O_{1}I$ of the triangle $A_{1}B_{1}C_{1}$ is their common radical axis. For any reference triangle $ABC$, its image $A_{1}B_{1}C_{1}$ under $\Psi$ is always acute. The chords $A_{1}J_{a}$, $B_{1}J_{b}$, $C_{1}J_{c}$ of the circles $\mathcal P_{a}$, $\mathcal P_{b}$, $\mathcal P_{c}$ meet at $I$ in their interiors, hence, these coaxal circles also intersect at two points $L_{1},$ $L_{1}^{\prime}$ on their common radical axis $O_{1}I$. Since $A_{1}O_{a}$, $B_{1}O_{b}$, $C_{1}O_{c}$ are diameters of the circles $\mathcal P_{a}$, $\mathcal P_{b}$, $\mathcal P_{c}$, the angles $\angle A_{1}L_{1}O_{a}$, $\angle B_{1}L_{1}O_{b}$, $\angle C_{1}L_{1}O_{c}$ are right and the line $O_{a}O_{b}O_{c}$ is the orthotransversal of the triangle $A_{1}B_{1}C_{1}$ with respect to the point $L_{1}$ or $L_{1}^{\prime}$ on its Euler line $O_{1}I$.

Let $K_{1}$ be the symmedian point of the inverted triangle $A_{1}B_{1}C_{1}$ and $K_{a}$, $K_{b}$, $K_{c}$ the feet of the symmedians through $A_{1}$, $B_{1}$, $C_{1}$. Because of the harmonic cross ratios

$$(O_{a}, K_{a}, B_{1}, C_{1}) = (O_{b}, K_{b}, C_{1}, A_{1}) = (O_{c}, K_{c}, A_{1}, B_{1}) = -1,$$

the line $O_{a}O_{b}O_{c}$ is the tripolar with respect to symmedian point $K_{1}$. Consequently, $K_{1}$ is the orthocorrespondent of the points $L_{1}$ and $L_{1}^{\prime}$, $L_{1}^{\perp} = L_{1}^{\prime\perp} = K_{1}$ [2].

\bigskip

{\bf{2. 4. The incenter}}

\bigskip

When the point $Q$ in Theorem 5 is identical with the incenter $I = I^{*}$, its own isogonal conjugate, the circles $(AXA^{\prime})$, $(BYB^{\prime})$, $(CZC^{\prime})$ degenerate into a pencil of concurrent lines, the internal bisectors $AX$, $BY$, $CZ$ of the angles $\angle A$, $\angle B$, $\angle C$. Thus, Theorem 5 is trivial for the incenter $I = I^{*}$. $\square$

The internal angle bisectors $AX$, $BY$, $CZ$ pass through the inversion center $I$, hence, the inversion $\Psi$ carries them into themselves. They become the altitudes $A_{1}X_{1}$, $B_{1}Y_{1}$, $C_{1}Z_{1}$ of the image $A_{1}B_{1}C_{1}$ of the reference triangle $ABC$ and the inversion center $I$ becomes its orthocenter. These altitudes can be considered as circles of infinite radius centered on the sidelines of the triangle $A_{1}B_{1}C_{1}$ and passing through its vertices $A_{1}$, $B_{1}$, $C_{1}$, respectively. The line at infinity is the orthotransversal with respect to the orthocenter $I$, which lies on the Euler line $O_{1}I$ of the triangle $A_{1}B_{1}C_{1}$. A pencil of concurrent lines does not have a radical axis; any line of the pencil would apply.

\bigskip

{\bf{3. Inversion in the incircle}}

\bigskip

From our analysis of the four special cases $X_{56} = M^{*}$, $X_{58} = S^{*}$, $K = G^{*}$ and $I = I^{*}$ with the help of the inversion in the triangle incircle, we conjecture the next theorem:

\bigskip
\begin{center}
\includegraphics[scale=0.8]{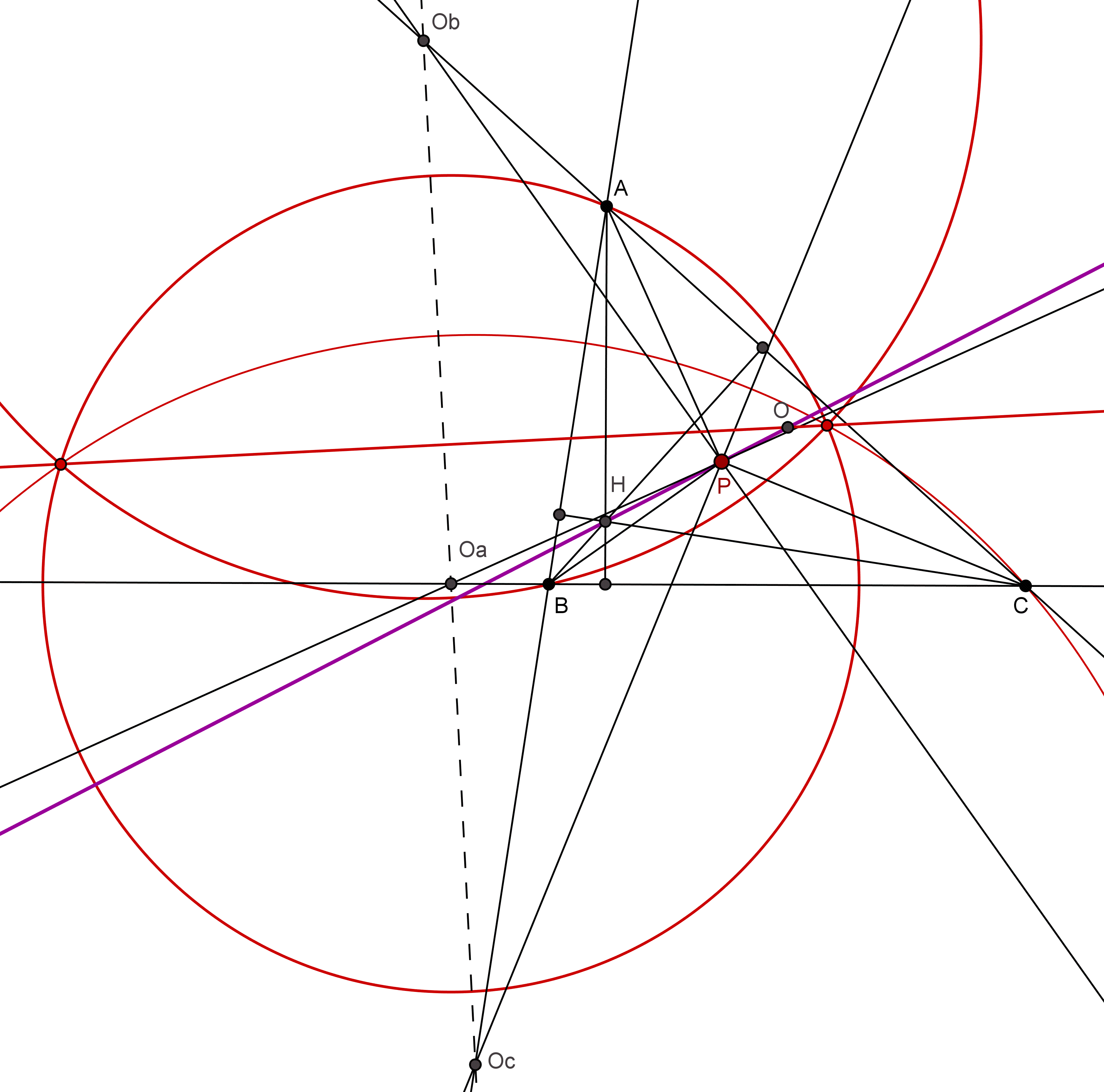}
\end{center}
\bigskip

{\bf{Theorem 7}}. Let $P$ be an arbitrary point in the plane of a non-equilateral triangle $ABC.$  Let $O_{a}O_{b}O_{c}$ be the orthotransversal of the point $P$, with the points $O_{a}$, $O_{b}$, $O_{b}$ lying on the triangle sidelines $BC$, $CA$, $AB$, respectively. The circles $\mathcal O_{a}$, $\mathcal O_{b}$, $\mathcal O_{c}$, centered at $O_{a}$, $O_{b}$, $O_{c}$ and passing through the triangle vertices $A$, $B$, $C$, respectively, are coaxal and the triangle circumcenter $O$ lies on their radical axis, if and only if the point $P$ lies on the Euler line of triangle $ABC$.

\bigskip

Before approaching the proof of Theorem 7, we need the following two simple lemmas. The second one is the proposition of existence of an orthotransversal and we present another proof of this well known fact [2]. 

\bigskip

{\bf{Lemma 8}}. Let $P$ be an arbitrary point in the plane of a triangle $ABC$. Let $A^{\prime}$, $B^{\prime}$, $C^{\prime}$ be the midpoints of the triangle sides $BC$, $CA$, $AB$. Let $\mathcal P_{a}$, $\mathcal P_{b}$, $\mathcal P_{c}$ be three circles centered on the triangle midlines $B^{\prime}C^{\prime}$, $C^{\prime}A^{\prime}$, $A^{\prime}B^{\prime}$, passing through the triangle vertices $A$, $B$, $C$ and intersecting at the point $P$. The circles $\mathcal P_{a}$, $\mathcal P_{b}$, $\mathcal P_{c}$ are coaxal and the triangle orthocenter $H$ lies on their common radical axis.

\bigskip

{\it{Proof of Lemma 8}}. Let $X$, $Y$, $Z$ be the feet of the triangle altitudes from the vertices $A$, $B$, $C$. Let the altitudes $AH$, $BH$, $CH$ meet the triangle circumcircle $\mathcal O$ again at points $H_{a}$, $H_{b}$, $H_{c}$. Power of the orthocenter $H$ to the circumcircle $\mathcal O$ is

$$\overline{HH_{a}} \cdot \overline{HA} = \overline{HH_{b}} \cdot \overline{HB} = \overline{HH_{c}} \cdot \overline{HC}.$$

Since the circles $\mathcal P_{a}$, $\mathcal P_{b}$, $\mathcal P_{c}$ are centered on the triangle midlines $B^{\prime}C^{\prime}$, $C^{\prime}A^{\prime}$, $A^{\prime}B^{\prime}$ and pass through the corresponding vertices $A$, $B$, $C$, they also pass through the corresponding altitude feet $X$, $Y$, $Z$. The altitude intersections $H_{a}$, $H_{b}$, $H_{c}$ with the circumcircle $\mathcal O$ other that the vertices $A$, $B$, $C$ are reflections of the orthocenter $H$ in the triangle sidelines $BC$, $CA$, $AB$. Powers of $H$ to the circles $\mathcal P_{a}$, $\mathcal P_{b}$, $\mathcal P_{c}$ are then

$$\overline{HX} \cdot \overline{HA} = \frac{1}{2} \overline{HH_{a}} \cdot \overline{HA},\ \overline{HY} \cdot \overline{HB} = \frac{1}{2} \overline{HH_{b}} \cdot \overline{HB},\ \overline{HZ} \cdot \overline{HC} = \frac{1}{2} \overline{HH_{c}} \cdot \overline{HC},$$

respectively. Thus, powers of the orthocenter $H$ to all three circles are the same, equal to half the power of $H$ to the circumcircle $\mathcal O$. But, these three circles intersect at the point $P$ and powers of $P$ to all three circles are zero, {\it{i.e.}}, also the same. Consequently, the line $PH$ is the common radical axis of the circles $\mathcal P_{a}$, $\mathcal P_{b}$, $\mathcal P_{c}$ and they are coaxal. $\square$
 
\bigskip

{\bf{Lemma 9}}. Let $O_{a}$, $O_{b}$, $O_{c}$ be the intersections of the circles $\mathcal P_{a}$, $\mathcal P_{b}$, $\mathcal P_{c}$, defined in Lemma 8, with the triangle sidelines $BC$, $CA$, $AB$ other than the altitude feet $X$, $Y$, $Z$. The points $O_{a}$, $O_{b}$, $O_{c}$ are collinear and the line $O_{a}O_{b}O_{c}$ is the orthotransversal of the triangle $ABC$ with respect to the point $P$.

\bigskip

{\it{Proof of Lemma 9}}. Since the angles $\angle AXO_{a}$, $\angle BYO_{b}$, $\angle CZO_{c}$ are right, the segments $AO_{a}$, $BO_{b}$, $CO_{c}$ are diameters of the circles $\mathcal P_{a}$, $\mathcal P_{b}$, $\mathcal P_{c}$. Consequently, the angles $\angle APO_{a}$, $\angle BPO_{b}$, $\angle CPO_{c}$ are also right. The circles $\mathcal P_{a}$, $\mathcal P_{b}$, $\mathcal P_{c}$ are centered on the midlines $B^{\prime}C^{\prime}$, $ C^{\prime}A^{\prime}$, $A^{\prime}B^{\prime}$ of the triangle $ABC$, identical with the sidelines of its medial triangle $A^{\prime}B^{\prime}C^{\prime}$. By Lemma 8, these three circles are coaxal, hence, their centers $P_{a}$, $P_{b}$, $P_{c}$ are collinear. Since $P_{a}$, $B^{\prime}$, $C^{\prime}$ are the midpoints of the segments $AO_{a}$, $CA$, $AB$, after cyclic exchange, it follows that

$$\frac{\overline{O_{a}B}}{\overline{O_{a}C}} = \frac{\overline{P_{a}C^{\prime}}}{\overline{P_{a}B^{\prime}}},\ \frac{\overline{O_{b}C}}{\overline{O_{b}A}} = \frac{\overline{P_{b}A^{\prime}}}{\overline{P_{b}C^{\prime}}},\ \frac{\overline{O_{c}A}}{\overline{O_{c}B}} = \frac{\overline{P_{c}B^{\prime}}}{\overline{P_{c}A^{\prime}}}.$$

By Menelaus theorem,

$$\frac{\overline{O_{a}B}}{\overline{O_{a}C}} \cdot \frac{\overline{O_{b}C}}{\overline{O_{b}A}} \cdot \frac{\overline{O_{c}A}}{\overline{O_{c}B}} = \frac{\overline{P_{a}C^{\prime}}}{\overline{P_{a}B^{\prime}}} \cdot \frac{\overline{P_{b}A^{\prime}}}{\overline{P_{b}C^{\prime}}} \cdot \frac{\overline{P_{c}B^{\prime}}}{\overline{P_{c}A^{\prime}}} = 1$$

and the points $O_{a}$, $O_{b}$, $O_{c}$ are also collinear. Because of the right angles $\angle APO_{a}$, $\angle BPO_{b}$, $\angle CPO_{c}$, the line $O_{a}O_{b}O_{c}$ is the orthotransversal of the triangle $ABC$ with respect to the point $P$. $\square$

\bigskip

{\it{Proof of Theorem 7 (sufficiency)}}. Let $\mathcal P_{a}$, $\mathcal P_{b}$, $\mathcal P_{c}$ be the circles with diameters $AO_{a}$, $BO_{b}$, $CO_{c}$. Since the line $O_{a}O_{b}O_{c}$ is the orthotransversal with respect to the point $P$, the angles $\angle APO_{a}$, $\angle BPO_{b}$, $\angle CPO_{c}$ are right and the circles $\mathcal P_{a}$, $\mathcal P_{b}$, $\mathcal P_{c}$ concur at $P$. By Lemma 8, they are coaxal with the common radical axis $PH$. By assumption, the point $P$ lies on the triangle Euler line $OH$. Thus, the circumcenter $O$ is on the common radical axis of these three circles and power of $O$ to all three circles is the same. Let the lines $AO$, $BO$, $CO$ meet the circles $\mathcal P_{a}$, $\mathcal P_{b}$, $\mathcal P_{c}$ again at points $U$, $V$, $W$. Since

$$\overline{OU} \cdot \overline{OA} = \overline{OV} \cdot \overline{OB} = \overline{OW} \cdot \overline{OC},$$

and $\overline{OA} = \overline{OB} = \overline{OC} = R$, where $R$ is the circumradius length, it follows that $\overline{OU} = \overline{OV} = \overline{OW}$. Consequently, the lines $VW$, $WU$, $UV$ are parallel to the respective triangle sidelines $BC$, $CA$, $AB$. Let the lines $AO$, $BO$, $CO$ meet the circles $\mathcal O_{a}$, $\mathcal O_{b}$, $\mathcal O_{c}$ again at points $X$, $Y$, $Z$. Since the circles $\mathcal O_{a}$, $\mathcal O_{b}$, $\mathcal O_{c}$ have radii $AO_{a}$, $BO_{b}$, $CO_{c}$ equal to diameters of the circles $\mathcal P_{a}$, $\mathcal P_{b}$, $\mathcal P_{c}$, these circle pairs are tangent at $A$, $B$, $C$. The circle pairs $\mathcal P_{a}$, $\mathcal O_{a}$; $\mathcal P_{b}$, $\mathcal O_{b}$, $\mathcal P_{c}$, $\mathcal O_{c}$ are then centrally similar with similarity centers $A$, $B$, $C$ and coefficient $\frac{1}{2}$, therefore, $U$, $V$, $W$ are the midpoints of the segments $AX$, $BY$, $CZ$, respectively. Thus, we conclude that the lines $YZ$, $ZX$, $XY$ are also parallel to the respective triangle sidelines $BC$, $CA$, $AB$ and $\overline{OX} = \overline{OY} = \overline{OZ}$. Powers of the circumcenter $O$ to the circles $\mathcal O_{a}$, $\mathcal O_{b}$, $\mathcal O_{c}$ are then equal,

$$\overline{OX} \cdot \overline{OA} = \overline{OY} \cdot \overline{OB} = \overline{OZ} \cdot \overline{OC}.$$

Since their centers $O_{a}$, $O_{b}$, $O_{c}$ are collinear, their pairwise radical axes are parallel. Existence of a point with the same power to all three circles then implies a common radical axis, {\it{i.e.}}, the circles $\mathcal O_{a}$, $\mathcal O_{b},$ $ \mathcal O_{c}$ are coaxal and the circumcenter $O$ lies on their common radical axis. $\square$

\bigskip

{\it{Proof of Theorem 7 (necessity)}}. Assume now that the point $P$ is not on the Euler line $OH$ and that the circles $\mathcal O_{a}$, $\mathcal O_{b}$, $\mathcal O_{c}$ centered at the intersections $O_{a}$, $O_{b}$, $O_{c}$ of the triangle sidelines with the orthotransversal with respect to the point $P$ are still coaxal. Let the perpendicular bisectors of the triangle sides $BC$, $CA$, $AB$ cut the radical axis $PH$ of the circles $\mathcal P_{a}$, $\mathcal P_{b}$, $\mathcal P_{c}$ with the diameters $AO_{a}$, $BO_{b}$, $CO_{c}$ at points $R_{a}$, $R_{b}$, $R_{c}$. Exactly as in the proof of sufficiency of Theorem 7, we can show that the points $R_{a}$, $R_{b}$, $R_{c}$ are on the radical axes of the circle pairs $\mathcal O_{b}$, $\mathcal O_{c}$; $\mathcal O_{c}$, $ \mathcal O_{a}$; $\mathcal O_{a}$, $\mathcal O_{b}$. Since these three circles are assumed to be coaxal, the line $R_{a}R_{b}R_{c}$, identical with the line $PH$, is their common radical axis and powers of the orthocenter $H$ to all three circles are equal. But, the circles $\mathcal O_{a}$, $\mathcal O_{b}$, $\mathcal O_{c}$, centered on the triangle sidelines and passing through the vertices $A$, $B$, $C$, also pass through their reflections $A^{\prime\prime}$, $B^{\prime\prime}$, $C^{\prime\prime}$ in the opposite sidelines. Consequently, powers of $H$ to $\mathcal O_{a}$, $\mathcal O_{b}$, $\mathcal O_{c}$ are equal to

$$\overline{HA} \cdot \overline{HA^{\prime\prime}} = \overline{HB} \cdot \overline{HB^{\prime\prime}} = \overline{HC} \cdot \overline{HC^{\prime\prime}},$$

$$\overline{HA} \cdot (\overline{HH_{a}} + \overline{H_{a}A^{\prime\prime}}) = \overline{HB} \cdot  (\overline{HH_{b}} + \overline{H_{b}B^{\prime\prime}}) = \overline{HC} \cdot (\overline{HH_{c}} + \overline{H_{c}C^{\prime\prime}}),$$

where $H_{a}$, $H_{b}$, $H_{c}$ are the altitude intersections with the circumcircle $\mathcal O$ other that the vertices $A$, $B$, $C$. Subtracting the power of $H$ to $\mathcal O$ and substituting $\overline{H_{a}A^{\prime\prime}} = -\overline{HA}$, $\overline{H_{b}B^{\prime\prime}} = -\overline{HB}$, $\overline{H_{c}C^{\prime\prime}} = -\overline{HC}$ leads to

$$\overline{HA}^{2} = \overline{HB}^{2} = \overline{HC}^2.$$
 
Solutions, such as $\overline{HA} = \overline{HB} = -\overline{HC}$, leading to $\angle A = \angle B = 0$, $\angle C = \pi$, are not acceptable. This implies that all triangle altitudes are equal and the triangle $ABC$ is equilateral, which is a contradiction. $\hfill \square$
 
\bigskip

{\it{Remark}}. For an equilateral triangle, the Euler line is undefined and Theorem 7 fails. Since any line through the equilateral triangle circumcenter can be considered as its Euler line, the circles $\mathcal O_{a}$, $\mathcal O_{b}$, $\mathcal O_{c}$, centered on the orthotransversal $O_{a}O_{b}O_{c}$ with respect to any point $P$ in the triangle plane, are always coaxal and their common radical axis passes through the equilateral triangle circumcenter.

\bigskip

In section 2, we have demonstrated that for four points $X_{56}$, $X_{58}$, $K$ and $I$ lying of the circumhyperbola $\mathcal H = (IG)^{*}$ of the reference triangle $ABC$, the inversion $\Psi$ with center $I$ and power $r^2$ carried the coaxality claim of Theorem 5 into the coaxality claim of Theorem 7 for the image $A_{1}B_{1}C_{1}$ of the reference triangle $ABC$ under $\Psi$. Now we use Theorem 7 to prove the coaxality claim of Theorem 5 in general.

\bigskip

{\it{Proof of theorem 5 - coaxality}}. Let $XYZ$ be the cevian triangle of the incenter $I$ of the reference triangle $ABC$ and let $A^{\prime}B^{\prime}C^{\prime}$ be the circumcevian triangle of a point $Q$ lying of the circumhyperbola $\mathcal H = (IG)^{*}$. The circles $(AXA^{\prime})$, $(BYB^{\prime})$, $(CZC^{\prime})$ cut the triangle sidelines $BC$, $CA$, $AB$ again at points $U$, $V$, $W$. By Lemma 3, the feet $A_{0}$, $B_{0}$, $C_{0}$ of the perpendiculars from $U$, $V$, $W$ to the internal bisectors $AX$, $BY$, $CZ$ of the angles $\angle A$, $\angle B$, $\angle C$ form a triangle $A_{0}B_{0}C_{0}$, centrally similar to the reference triangle $ABC$ with similarity center $I$. Hereby,

$$\frac{\overline{AA_{0}}}{\overline{AI}} = \frac{\overline{BB_{0}}}{\overline{BI}} = \frac{\overline{CC_{0}}}{\overline{CI}}.$$

Let $I_{a}$, $I_{b}$, $I_{c}$ be the excenters in the angles $\angle A$, $\angle B$, $\angle C$ of the reference triangle  $ABC$ and $r_{a}$, $r_{b}$, $r_{c}$ the corresponding exradii lengths. Let $\mathcal Q_{a}$, $\mathcal Q_{b}$, $\mathcal Q_{c}$ be circles tangent to the circles $(AXA^{\prime})$, $(BYB^{\prime})$, $(CZC^{\prime})$ at $A$, $B$, $C$ and passing through the excenters $I_{a}$, $I_{b}$, $I_{c}$, respectively. Let $\mathcal O$, $\mathcal O_{0}$ be the circumcircles of the reference triangle $ABC$ and of the excentral triangle $I_{a}I_{b}I_{c}$, with circumcenters $O$, $O_{0}$ and circumradii $R$, $R_{0}$. The reference triangle $ABC$ is the orthic triangle of the excentral triangle $I_{a}I_{b}I_{c}$, hence, $I$ is the orthocenter, $\mathcal O$ is the nine-point circle of $I_{a}I_{b}I_{c}$ and $O_{0}OI$ its Euler line. Since the circles $\mathcal Q_{a}$, $\mathcal Q_{b}$, $\mathcal Q_{c}$ pass through the vertices $I_{a},$ $I_{b}$, $I_{c}$ and altitude feet $A$, $B$, $C$ of the excentral triangle, powers of its orthocenter $I$ to these three circles are equal to half the power of $I$ to the circumcircle $\mathcal O_{0}$, 

$$\overline{IA} \cdot \overline{II_{a}} = \overline{IB} \cdot \overline{II_{b}} = \overline{IC} \cdot \overline{II_{c}}.$$

Both the tangent of the circumcircle $\mathcal O_{0}$ of the excentral triangle $I_{a}I_{b}I_{b}$ at the vertex $I_{a}$ and the sideline $BC$ of the reference triangle $ABC$ are antiparallel to the line $I_{b}I_{c}$ with respect to the angle $\angle I_{c}I_{a}I_{b}$, hence, these two lines are parallel. Similarly, the tangents to $\mathcal O_{0}$ at $I_{b}$, $I_{c}$ are parallel to $CA$, $AB$. Let these tangents cut the circles $\mathcal Q_{a}$, $\mathcal Q_{b}$, $\mathcal Q_{c}$ again at points $U^{\prime}$, $V^{\prime}$, $W^{\prime}$. Since the lines $XU$, $I_{a}U^{\prime}$ are parallel, the line $UU^{\prime}$ passes through the tangency point $A$, the similarity center of the circles $(AXA^{\prime})$, $\mathcal Q_{a}$. 

Let the sidelines $I_{b}I_{c}$, $I_{c}I_{a}$, $I_{a}I_{b}$ of the excentral triangle $I_{a}I_{b}I_{c}$ cut the circles $(AXA^{\prime})$, $(BYB^{\prime})$, $(CZC^{\prime})$ at points $R$, $S$, $T$ and the circles $\mathcal Q_{a}$, $\mathcal Q_{b}$, $\mathcal Q_{c}$ at points $R^{\prime}$, $S^{\prime}$, $T^{\prime}$. The triangles $XUR$, $I_{a}U^{\prime}R^{\prime}$ are then centrally similar with similarity center $A$. Using the harmonic cross ratio $(A, X, I, I_{a}) = -1$ and the central similarity of the incircle $\mathcal I$ and the excircle $\mathcal I_{a}$ with similarity center $A$, their similarity coefficient is

$$\frac{\overline{AX}}{\overline{AI_{a}}} = \frac{2 \overline{AI}}{\overline{AI} + \overline{AI_{a}}} = \frac{2r}{r + r_{a}}.$$

In addition, the angle $\angle XAR = \angle I_{a}AR$ is right, hence, $XR$, $I_{a}R^{\prime}$ are diameters of the circles $(AXA^{\prime})$, $\mathcal Q_{a}$ and the angles $\angle XUR$, $\angle I_{a}U^{\prime}R^{\prime}$ are also right. By cyclic exchange, the right angle triangles $YVS$, $I_{b}V^{\prime}S^{\prime}$ are centrally similar with similarity center $B$ and coefficient $\frac{2r}{r + r_{b}}$ and the right angle triangles $ZWT$, $I_{c}W^{\prime}T^{\prime}$ are centrally similar with similarity center $C$ and coefficient $\frac{2r}{r + r_{c}}$. Consequently,

$$\frac{\overline{R^{\prime}U^{\prime}}}{\overline{S^{\prime}V^{\prime}}} = \frac{\overline{RU}}{\overline{SV}} \cdot \frac{r + r_{a}}{r + r_{b}}.$$

Consider the internal bisectors $AI_{a}$, $BI_{b}$, $CI_{c}$ of the angles $\angle A$, $\angle B$, $\angle C$ of the reference triangle $ABC$ as the corresponding altitudes of the excentral triangle $I_{a}I_{b}I_{c}$ with the orthocenter $I$. These altitudes meet the circumcircle $\mathcal O_{0}$ of the excentral triangle at the reflections $K_{a}$, $K_{b}$, $K_{c}$ of the orthocenter $I$ in the excentral triangle sidelines. Thus, we have

$$\frac{\overline{RU}}{\overline{AA_{0}}} = \frac{2R_{0}}{\overline{K_{a}I_{a}}} = \frac{2R_{0}}{\overline{AI} + \overline{AI_{a}}}$$

and by cyclic exchange,

$$\frac{\overline{SV}}{\overline{BB_{0}}} = \frac{2R_{0}}{\overline{BI} + \overline{BI_{b}}},\ \frac{\overline{TW}}{\overline{CC_{0}}} = \frac{2R_{0}}{\overline{CI} + \overline{CI_{c}}}.$$

Using the central similarity of the triangles $A_{0}B_{0}C_{0}$, $ABC$ with similarity center $I$ and the central similarity of the incircle $\mathcal I$ and the excircles $\mathcal I_{a}$, $\mathcal I_{b}$ with similarity centers $A$, $B$, we obtain

$$\frac{\overline{RU}}{\overline{SV}} = \frac{\overline{AI}}{\overline{BI}} \cdot \frac{\overline{BI} + \overline{BI_{b}}}{\overline{AI} + \overline{AI_{a}}} = \frac{r + r_{b}}{r + r_{a}}$$
 
With cyclic exchange, this yields $\overline{R^{\prime}U^{\prime}} = \overline{S^{\prime}V^{\prime}} = \overline{T^{\prime}W^{\prime}}$. Let the radii $O_{0}I_{a}$, $O_{0}I_{b}$, $O_{0}I_{c}$ of the circumcircle $\mathcal O_{0}$ of the excentral triangle $I_{a}I_{b}I_{c}$ meet the circles $\mathcal Q_{a}$, $\mathcal Q_{b}$, $\mathcal Q_{c}$ again at points $L_{a},$ $L_{b}$, $L_{c}$. The quadrilaterals $I_{a}L_{a}R^{\prime}U^{\prime}$, $I_{b}L_{b}S^{\prime}V^{\prime}$, $I_{c}L_{c}T^{\prime}W^{\prime}$ are all rectangles and consequently, $\overline{I_{a}L_{a}} = \overline{I_{b}L_{b}} = \overline{I_{c}L_{c}}$. As a result, powers of the circumcenter $O_{0}$ to the circles $\mathcal Q_{a}$, $\mathcal Q_{b}$, $\mathcal Q_{c}$ are equal, 

$$\overline{O_{0}L_{a}} \cdot \overline{O_{0}I_{a}} = \overline{O_{0}L_{b}} \cdot \overline{O_{0}I_{b}} = \overline{O_{0}L_{c}} \cdot \overline{O_{0}I_{c}}.$$

Since powers of both the orthocenter $I$ and circumcenter $O_{0}$ of the excentral triangle $I_{a}I_{b}I_{c}$ to the circles $\mathcal Q_{a}$, $\mathcal Q_{b}$, $\mathcal Q_{c}$ are equal, these three circles are coaxal and the Euler line $O_{0}I$ of the triangle $I_{a}I_{b}I_{c}$, identical with the diacentral line $OI$ of the reference triangle $ABC$, is their common radical axis. The chords $AI_{a}$, $BI_{b}$, $CI_{c}$ of the circles $\mathcal Q_{a}$, $\mathcal Q_{b}$, $\mathcal Q_{c}$ meet at their interiors at $I$, hence, these coaxal circles also intersect at two points $L$, $L^{\prime}$ on their common radical axis $OI$.

The inversion $\Psi$ with center $I$ and power $r^{2}$ carries the line $OI$ passing through the inversion center $I$ into itself. The excenters $I_{a}$, $I_{b}$, $I_{c}$ lie on the circumcircles $\mathcal M$, $\mathcal N$, $\mathcal P$ of the triangles $IBC$, $ICA$, $IAB$. The images of the circles $\mathcal M$, $\mathcal N$, $\mathcal P$ under $\Psi$ are the lines $B_{1}C_{1}$, $C_{1}A_{1}$, $A_{1}B_{1}$ and the images $J_{a}$ $J_{b}$, $J_{c}$ of the excenters $I_{a}$ $I_{b}$, $I_{c}$ are the respective altitude feet of the inverted triangle $A_{1}B_{1}C_{1}$. Consequently, the altitudes $A_{1}J_{a}$, $B_{1}J_{b}$, $C_{1}J_{c}$ are chords of the images $\mathcal P_{a}$, $\mathcal P_{b}$, $\mathcal P_{c}$ of the coaxal circles $\mathcal Q_{a}$, $\mathcal Q_{b}$, $\mathcal Q_{c}$. Thus, the coaxal circles $\mathcal P_{a}$, $\mathcal P_{b}$, $\mathcal P_{c}$ are centered on the corresponding midlines of the inverted triangle $A_{1}B_{1}C_{1}$ and they intersect on its Euler line $O_{1}I$, at the images $L_{1}$, $L_{1}^{\prime}$ of the points $L$, $L^{\prime}$. These circles cut the inverted triangle sidelines $B_{1}C_{1}$, $C_{1}A_{1}$, $A_{1}B_{1}$ again at points $O_{a}$, $O_{b}$, $O_{c}$. Since the angles $\angle A_{1}J_{a}O_{a}$, $\angle B_{1}J_{b}O_{b}$, $\angle C_{1}J_{c}O_{c}$ are right, the segments $A_{1}O_{a}$, $B_{1}O_{b}$, $C_{1}O_{c}$ are diameters of the circles $\mathcal P_{a}$, $\mathcal P_{b}$, $\mathcal P_{c}$. The angles $\angle A_{1}L_{1}O_{a}$, $\angle C_{1}L_{1}O_{b}$, $\angle C_{1}L_{1}O_{c}$ and $\angle A_{1}L_{1}^{\prime}O_{a}$, $\angle C_{1}L_{1}^{\prime}O_{b}$, $\angle C_{1}L_{1}^{\prime}O_{c}$ are then also right. By Lemma 9, the line $O_{a}O_{b}O_{c}$ is the common orthotransversal of the points $L_{1}$, $L_{1}^{\prime}$.
 
The images $X_{1}$, $Y_{1}$, $Z_{1}$ of the internal angle bisector feet $X$, $Y$, $Z$ of the reference triangle $ABC$ under $\Psi$ are the reflections of the vertices $A_{1}$, $B_{1}$, $C_{1}$ of the inverted triangle $A_{1}B_{1}C_{1}$ in its respective sidelines. Since $A_{1}X_{1}$, $B_{1}Y_{1}$, $C_{1}Z_{1}$ are chords of the images $(A_{1}X_{1}A_{1}^{\prime})$, $(B_{1}Y_{1}B_{1}^{\prime})$, $(C_{1}Z_{1}C_{1}^{\prime})$ of the circles $(AXA^{\prime})$, $(BYB^{\prime})$, $(CZC^{\prime})$, the inverted circles are centered on the respective sidelines of the inverted triangle $A_{1}B_{1}C_{1}$. Since the circles $(AXA^{\prime})$, $(BYB^{\prime})$, $(CZC^{\prime})$ are tangent to the circles $\mathcal Q_{a}$, $\mathcal Q_{b}$, $\mathcal Q_{c}$ at the vertices $A$, $B$, $C$ of the reference triangle $ABC$, their images $(A_{1}X_{1}A_{1}^{\prime})$, $(B_{1}Y_{1}B_{1}^{\prime})$, $(C_{1}Z_{1}C_{1}^{\prime})$ are tangent to the inverted circles $\mathcal P_{a}$, $\mathcal P_{b}$, $\mathcal P_{c}$ at the vertices $A_{1}$, $B_{1}$, $C_{1}$ of the inverted triangle $A_{1}B_{1}C_{1}$. It follows that the inverted circles $(A_{1}X_{1}A_{1}^{\prime})$, $(B_{1}Y_{1}B_{1}^{\prime})$, $(C_{1}Z_{1}C_{1}^{\prime})$ are centered at the intersections $O_{a}$, $O_{b}$, $O_{c}$ of the orthotransversal $O_{a}O_{b}O_{c}$ with respect to the point $L_{1}$ on the Euler line $O_{1}I$ with the triangle sidelines $B_{1}C_{1}$, $C_{1}A_{1}$, $A_{1}B_{1}$. By Theorem 7, these three circles are coaxal. But, inversion in any circle carries a pencil of circles or lines into a pencil of circles or lines. As a result, the circles $(AXA^{\prime})$, $(BYB^{\prime})$, $(CZC^{\prime})$ are also coaxal. 

Conversely, if the circles $(AXA^{\prime})$, $(BYB^{\prime})$, $(CZC^{\prime})$ are coaxal, their images under $\Psi$ are also coaxal and by Theorem 7, they are centered on the orthotransversal $O_{a}O_{b}O_{c}$ with respect to a point $L_{1}$ on the Euler line $O_{1}I$ of the inverted triangle $A_{1}B_{1}C_{1}$. Their tangent circles $\mathcal P_{a}$, $\mathcal P_{b}$, $\mathcal P_{c}$ with diameters $A_{1}O_{a}$, $B_{1}O_{b}$, $C_{1}O_{c}$ then intersect at the point $L_{1}$. The inversion $\Psi$ carries the circles $\mathcal P_{a}$, $\mathcal P_{b}$, $\mathcal P_{c}$ back to the coaxal circles $\mathcal Q_{a}$, $\mathcal Q_{b}$, $\mathcal Q_{c}$ tangent to the circles $(AXA^{\prime})$, $(BYB^{\prime})$, $(CZC^{\prime})$ at $A$, $B$, $C$, passing through the excenters $I_{a}$, $I_{b}$, $I_{c}$ and intesecting on the Euler line $O_{0}I$ of the excentral triangle $I_{a}I_{b}I_{c}$, identical with the diacentral line $OI$ of the reference triangle $ABC$. By reversing the remaining arguments of the direct proof, we arrive back to the triangle $A_{0}B_{0}C_{0}$ centrally similar the reference triangle $ABC$ with similarity center $I$, the external angle bisectors of which intersect the respective sidelines $BC$, $CA$, $AB$ of the reference triangle $ABC$ at points $U$, $V$, $W$, lying on the circles $(AXA^{\prime})$, $(BYB^{\prime})$, $(CZC^{\prime})$. By Lemma 3, the lines $AA^{\prime}$, $BB^{\prime}$, $CC^{\prime}$ then intersect at a point $Q$ on the circumhyperbola $\mathcal H = (IG)^{*}$.  $\hfill \square$

\bigskip

{\bf{4. General case}}

\bigskip

{\bf{Theorem 10}}. Let $P$, $Q$ be two arbitrary point in the plane of a non-equilateral triangle $ABC$, different from the triangle vertices, the first one also different from the triangle orthocenter, $P \neq H$. Let $XYZ$ be the cevian triangle of the point $P$ and $A^{\prime}B^{\prime}C^{\prime}$ the circumcevian triangle of the point $Q$. The circles $(AXA^{\prime})$, $(BYB^{\prime})$, $(CZC^{\prime})$ are coaxal and their common radical axis passes through the isogonal conjugate $R = (P_{C})^{*}$ of the complement $P_{C}$ of the point $P$, if and only if the point $Q$ lies on the triangle circumconic $\mathcal K$ passing through the points $P$, $R$, the isogonal conjugate of the line $P^{*}P_{C}$.

\bigskip
\begin{center}
\includegraphics[scale=0.9]{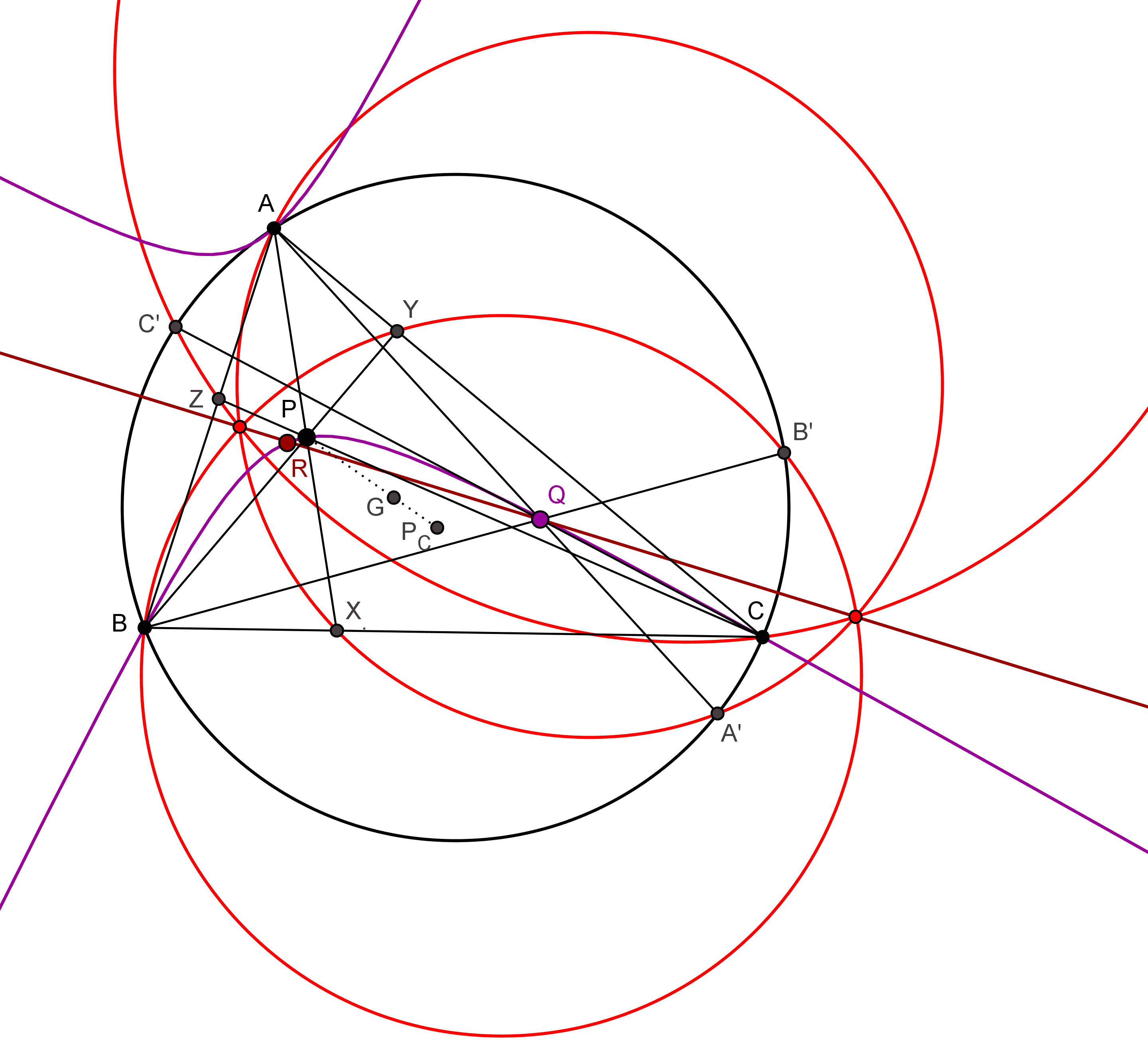}
\end{center}
\bigskip

{\it{Remark}}. Theorem 5 is a special case of Theorem 10 for $P = I$. Indeed, the incenter $I = I^{*}$ is its own isogonal conjugate, its complement is the Spiker point $S = I_{C}$ and isogonal conjugate of the Spiker point is the Kimberling center $X_{58} = S^{*} = (I_{C})^{*}$. The isogonal conjugate of the line $I^{*}I_{C} = IS = IG$ is then the circumhyperbola $\mathcal H = (IG)^{*}$. Theorem 5 fails, when the incenter and orthocenter are identical, $I = H$, {\it{i.e.}}, for an equilateral triangle.

\bigskip

{\it{Proof of Theorem 10}}. Let the points $P$, $Q$ have the barycentric coordinates $P = (p : q : r)$ and $Q = (u : v : w)$ and let $P^{*}$, $Q^{*}$ be their isogonal conjugates with respect to the triangle $ABC$. The equation of the line $P^{*}Q^{*}$ is then

$$\left(\frac{b^2c^2}{qw} - \frac{c^2b^2}{rv}\right) x + \left(\frac{c^2a^2}{ru} - \frac{a^2c^2}{pw}\right) y +\left(\frac{a^2b^2}{pv} - \frac{b^2a^2}{qu}\right) z = 0.$$

After multiplying by $(pqr) \cdot (uvw)$, the equation of the line $P^{*}Q^{*}$ can be written as

$$\left| \begin{array}{ccc}
a^{2}qr & b^{2}rp & c^{2}pq \\
a^{2}vw & b^{2}wu & c^{2}uv \\
x & y & z \\
\end{array}\right| = 0.$$

Assume that point $Q$ does not lie on the lines $AP$, $BP$, $CP$, {\it{i.e.}}, $wq - vr \neq 0$, $ur - wp \neq 0$, $vp - uq \neq 0$. The triangle circumcircle $\mathcal O$ has the equation $a^{2}yz + b^{2}zx + c^{2}xy = 0$. The radical axis $AQ$ of the circumcircle $\mathcal O$ and circle $(AXA^{\prime})$ has the equation $wy - vz = 0$. Hereby, the equation of the circle $(AXA^{\prime})$ is [8]

$$a^2yz + b^2zx + c^2xy + (wy - vz) (x + y + z) t = 0$$

for some $t$. Since the circle $(AXA^{\prime})$ passes through the trace $X = (0 : q : r)$ of the point $P$ in the triangle sideline $BC$,

$$t = -\frac{a^{2}qr}{(q + r)(wq - vr)}$$

and the equation of the circle $(AXA^{\prime})$ becomes

$$a^2yz + b^2zx + c^2xy - (x + y + z) \frac{a^{2}qr(wy - vz)}{(q + r)(wq - vr)} = 0.$$

By cyclic exchange, we get the equations of the circles $(BYB^{\prime})$, $(CZC^{\prime})$. Subtracting the equations of the circles $(BYB^{\prime})$, $(CZC^{\prime})$ from each other yields the equation of their radical axis $r_{BC}$:

$$r_{BC}(x, y, z) = \frac{c^{2}q(vx - uy)}{(p + q)(vp - uq)} - \frac{b^{2}r(uz - wx)}{(r + p)(ur - wp)} = 0.$$

Cyclically, we obtain the equations of the radical axes $r_{CA}$, $r_{AB}$ of the circle pairs $(CZC^{\prime})$, $(AXA^{\prime})$; $(AXA^{\prime})$, $(BYB^{\prime})$.

Powers of the point $Q$ to the circles  $(AXA^{\prime})$, $(BYB^{\prime})$, $(CZC^{\prime})$ are equal to the power of $Q$ to the circumcircle $\mathcal O$,

$$\overline{QA} \cdot \overline{QA^{\prime}} = \overline{QB} \cdot \overline{QB^{\prime}} = \overline{QC} \cdot \overline{QC^{\prime}},$$

hence, $Q$ is at least the radical center of the circles $(AXA^{\prime})$, $(BYB^{\prime})$, $(CZC^{\prime})$ and their pairwise radical axes $r_{BC}$, $r_{CA}$, $r_{AB}$ concur at $Q$.

Assume that the complement $P_{C}$ of $P$ lies on the line $P^{*}Q^{*}$. The point $Q$ then lies on the triangle circumconic $\mathcal K$ passing through the points $P$ and $R = (P_{C})^{*}$, the isogonal conjugate of the line $P^{*}P_{C}$. This is equivalent to

$$\Delta = \left| \begin{array}{ccc}
a^{2}qr & b^{2}rp & c^{2}pq \\
a^{2}vw & b^{2}wu & c^{2}uv \\
q + r & r + p & p + q \\
\end{array}\right| = 0.$$

Substituting the barycentric coordinates of the point

$$R = (P_{C})^{*} = \left(\frac{a^{2}}{q + r} : \frac{b^{2}}{r + p} : \frac{c^{2}}{p + q}\right)$$

into the equation of the radical axis $r_{BC}$, the left side becomes proportional to the determinant $\Delta$:

$$r_{BC}(R) = \frac{\Delta}{(q + r)(r + p)(p + q)(vp - uq)(ur - wp)}.$$

Since $\Delta = 0$, the point $R$ lies on the radical axis $r_{BC}$ and similarly, $R$ lies on the radical axes $r_{CA}$, $r_{AB}$. But, these radical axes concur at the point $Q$, generally different from the point $R$, hence, they are identical. The line $QR$ is then the common radical axis of the coaxal circles $(AXA^{\prime})$, $(BYB^{\prime})$, $(CZC^{\prime})$. If the points $Q$, $R$ are different from each other, the common radical axis intersects the circumconic $\mathcal K$ at $Q$, $R$. If the points $Q$, $R$ coincide, the common radical axis is tangent to the circumconic $\mathcal K$ at $R$ and the circles $(AXA^{\prime})$, $(BYB^{\prime})$, $(CZC^{\prime})$ remain coaxal by the continuity principle.

Conversely, assume that the circles $(AXA^{\prime})$, $(BYB^{\prime})$, $(CZC^{\prime})$ are coaxal, {\it{i. e.}}, their pairwise radical axes $r_{BC}$, $r_{CA}$, $r_{AB}$ are identical with their common radical axis. Any non-trivial linear combination of their equations is also an equation of the same line. But, substituting the barycentric coordinates of the point $R$ into the left side of a linear combination, such as

$$\frac{r_{BC}(x, y, z)}{wq - vr} - \frac{r_{CA}(x, y, z)}{ur - wp} = 0,$$  

yields zero, which means that the point $R$ lies on the common radical axis. This implies that the determinant $\Delta = 0$, the isogonal conjugate $R$ of the complement $P_{C}$ of the point $P$ lies on the line $P^{*}Q^{*}$ and the point $Q$ on the circumconic $\mathcal K$ through the points $P$, $R$.

If the line $P^{*}P_{C}$ passes through any triangle vertex, say $A$, the circumconic $\mathcal K$ degenerates into the line $APR$. By the continuity principle, the circles $(AXA^{\prime})$, $(BYB^{\prime})$, $(CZC^{\prime})$ are coaxal, if and only if the point $Q$ lies on the line $APR$. The circle $(AXA^{\prime})$ then degenerates into the line $AXA^{\prime}$, identical with the line $APR$ and with the radical axis of the circles $(BYB^{\prime})$, $(CZC^{\prime})$. $\hfill \square$
\bigskip

{\it{Remark}}. Both the isogonal conjugate $H^{*}$ and complement $H_{C}$ of the orthocenter $H$ are identical with the circumcenter $O$. The orthocenter is the only point, the isogonal conjugate and complement of which are identical. When $P = H$, the line $P^{*}P_{C}$ and the circumconic $\mathcal K$ are undefined. The complement $H_{C} = O$ of $H$ is on the line $H^{*}Q^{*} = OQ^{*}$ for any point $Q$, the circles $(AXA^{\prime})$, $(BYB^{\prime})$, $(CZC^{\prime})$ are always coaxal and their common radical axis is the line $QH$. Powers of the orthocenter $H$ to all three circles are equal to half the power of $H$ to the circumcircle $\mathcal O$.

\bigskip

\bigskip

\bigskip

\bigskip

\bigskip

\bigskip

{\bf{References}}

\bigskip

[1] N. A. Court, On the Circles of Apollonius, {\it{Amer. Math. Monthly}} {\bf{22}}, No. 8. (Oct. 1915), pp. 261-263.

\bigskip

[2] B. Gibert, Orthocorrespondence and orthopivotal cubics, {\it{Forum Geom.}} {\bf{3}} (2003), pp. 1-27. 

\bigskip

[3] A. P. Hatzipolakis, et al., Concurrency of Four Euler Lines, {\it{Forum Geom.}} {\bf{1}}(2001), pp. 59-68. 

\bigskip

[4] C. Kimberling, {\it{Encyclopedia of Triangle Centers}}, available at

\verb+http://faculty.evansville.edu/ck6/encyclopedia/ETC.html+.

\bigskip

[5] K. L. Nguyen and J. C. Salazar, On the mixtilinear incircles and excircles, {\it{Forum Geom.}} {\bf{6}} (2006), pp. 1-16. 

\bigskip

[6] C. Pohoata, P. Yiu, Hyacinthos messages 15606, 15607, October 6-7, 2007,

\verb+http://tech.groups.yahoo.com/group/Hyacinthos+

\bigskip

[7] P. Yiu, Mixtilinear Incircles, {\it{Amer. Math. Monthly}} {\bf{106}}, No. 10 (Dec. 1999), pp. 952-955.

\bigskip

[8] P. Yiu, Introduction to the Geometry of the Triangle, {\it{Florida Atlantic University Lecture Notes}}, 2001. 

\bigskip

\bigskip

\bigskip

Cosmin Pohoata: 13 Pridvorului Street, Bucharest 010014, Romania.

{\it{E-mail address}}: pohoata\_{c}osmin2000@yahoo.com

\bigskip

Vladimir Zajic: 29 Milldown Road, Yaphank, NY 11980, U.S.A.

Brookhaven National Laboratory, Upton, NY 11973, U.S.A. 

{\it{E-mail address}}: vzajic@bnl.gov

\end{document}